\documentclass [a4paper, 12pt]{amsart}

\usepackage{times}

\usepackage{amsmath,amsfonts,amsthm,amssymb,amscd}
\usepackage{graphicx, epsfig, psfrag}
\usepackage{pb-diagram}
\usepackage{hyperref}
\usepackage{wasysym, auto-pst-pdf}
\usepackage{epstopdf,yfonts}
\usepackage{fancyhdr}
\usepackage{tikz}
\usepackage[all]{xy}
\usepackage{color}
\usepackage{xargs}
\usepackage{comment}   
\usepackage[sc]{mathpazo}
\usepackage{eulervm}

\theoremstyle{plain}
\newtheorem{cor}{Corollary}[section]
\newtheorem{thm}[cor]{Theorem}
\newtheorem{prop}[cor]{Proposition}
\newtheorem{q}[cor]{Question}

\theoremstyle{definition}
\newtheorem{defi}[cor]{Definition}
\theoremstyle{remark}
\newtheorem{remark}[cor]{Remark}

\newtheorem{ex}[cor]{Example}

\newcommand\bqn{\begin{equation*}}
\newcommand\eqn{\end{equation*}}
\newcommand\bq{\begin{equation}}
\newcommand\eq{\end{equation}}
\newcommand\be{\begin{enumerate}}
\newcommand\ee{\end{enumerate}}
\newcommand\bei{\begin{itemize}}
\newcommand\eei{\end{itemize}}
\newcommand\ba{\begin{aligned}}
\newcommand\ea{\end{aligned}}
\newcommand\ban{\begin{aligned*}}
\newcommand\ean{\end{aligned*}}
\newcommand{\bsm}{\left(\begin{smallmatrix}}
\newcommand{\esm}{\end{smallmatrix}\right)}                   
\newcommand{\bpm}{\begin{pmatrix}}
\newcommand{\epm}{\end{pmatrix}}

\newcommand{\thismonth}{\ifcase\month 
  \or January\or February\or March\or April\or May\or June%
  \or July\or August\or September\or October\or November%
  \or December\fi}

\newcommand{\calA}{\mathcal A}
\newcommand{\calB}{\mathcal B}
\newcommand{\calC}{\mathcal C}

\newcommand{\calE}{\mathcal E}
\newcommand{\calF}{\mathcal F}

\newcommand{\calH}{\mathcal H}

\newcommand{\calL}{\mathcal L}

\newcommand{\calP}{\mathcal P}
\newcommand{\calQ}{\mathcal Q}
\newcommand{\calR}{\mathcal R}
\newcommand{\calS}{\mathcal S}
\newcommand{\calT}{\mathcal T}

\newcommand{\calV}{\mathcal V}

\newcommand{\calX}{\mathcal X}

\newcommand{\fa}{\mathfrak a}

\newcommand{\fp}{\mathfrak p}

\newcommand{\CC}{\mathbb C}

\newcommand{\FF}{\mathbb F}

\newcommand{\KK}{\mathbb K}
\newcommand{\LL}{\mathbb L}

\newcommand{\NN}{\mathbb N}

\newcommand{\PP}{\mathbb P}
\newcommand{\QQ}{\mathbb Q}
\newcommand{\RR}{\mathbb R}

\newcommand{\ZZ}{\mathbb Z}

\newcommand{\bA}{\mathbf A}

\newcommand{\bG}{\mathbf G}

\newcommand{\SL}{\operatorname{SL}}
\newcommand{\Sp}{\operatorname{Sp}}

\newcommand{\PSL}{\operatorname{PSL}}

\newcommand{\SO}{\operatorname{SO}}

\newcommand\Sym{\mathrm{Sym}}

\newcommand{\Out}{\operatorname{Out}}
\newcommand\T{\operatorname{T}}
\newcommand{\mcgp}{\mathcal{MCG}^+(\Sigma)}
\newcommand{\Iso}{\operatorname{Iso}}

\newcommand{\Hom}{\mathrm{Hom}}
\newcommand{\Homred}{{\mathrm{Hom}_\mathrm{red}}}
\newcommand{\Homhit}{{\mathrm{Hom}_\mathrm{Hit}}}

\newcommand{\Ad}{\mathrm{Ad}}
\newcommand\diag{\mathrm{diag}}
\newcommand\Id{\mathrm{Id}}
\newcommand\ov{\overline}
\newcommand\per{\mathrm{per}}
\newcommand\pis{{\pi_1(\Sigma)}}
\newcommand\pish{{\pi_1(\Sigma)_\mathrm{h}}}
\newcommand\Stab{\mathrm{Stab}}
\newcommand\tr{\mathrm{tr}}
\newcommand\wt{\widetilde}

\newcommand\dnf{d_{N_\mathbb F}}
\newcommand\Rspec{\operatorname{Spec}_\mathrm{r}}
\newcommand\Syst{\mathrm{Syst}}
\newcommand\ftn{\mathbb F^{2n}}
\newcommand\ktn{\mathbb K^{2n}}
\newcommand{\qbarr}{\overline{\mathbb Q}^\mathrm{r}}
\newcommand\rol{{\mathbb R^\omega_{\pmb{\lambda}}}}
\newcommand\rhol{{\rho^\omega_{\pmb{\lambda}}}}
\newcommand\ol{{\mathcal{O}^\omega_{\pmb{\lambda}}}}
\newcommand\kol{{k^\omega_{\pmb{\lambda}}}}
\newcommand\romu{{\mathbb R^\omega_{\pmb{\mu}}}}
\newcommand\rhomu{{\rho^\omega_{\pmb{\mu}}}}
\newcommand\Ximax{\Xi_\mathrm{max}}
\newcommand\abp{\overline{A}^+}
\newcommand\fabp{\overline{\mathfrak a}^+}
\newcommand\XiHit{\Xi_\mathrm{Hit}} 
\newcommand\bgf{\mathcal B_{G_\mathbb F}}
\newcommand\bsp{\mathcal B_{\Sp(2n,\mathbb F)}}
\newcommand\bpsl{\mathcal B_{\PSL(2n,\mathbb F)}}
\newcommand\obpsl{\overline{\mathcal B_{\PSL(2n,\mathbb F)}}}
\newcommand\bgrol{\mathcal B_{G_{{\mathbb R^\omega_{\pmb{\lambda}}}}}}
\newcommand\bgromu{\mathcal B_{G_{{\mathbb R^\omega_{\pmb{\mu}}}}}}
\newcommand\bspkol{\mathcal B_{\Sp(2n,{k^\omega_{\pmb{\lambda}}})}}
\newcommand\obgf{\overline{\mathcal B_{G_\mathbb F}}}
\newcommand\obsp{\overline{\mathcal B_{\Sp(2n,\mathbb F)}}}
\newcommand\kkb{\kappa_\mathbb K^\mathrm{b}}
\newcommand\hb{\mathrm{H}_\mathrm{b}}
\def\rsp#1{#1^\mathrm{RSp}}
\def\rspcl#1{#1^\mathrm{RSp}_\mathrm{cl}}
\newcommand\rspclmax{\Xi_\mathrm{max, cl}^\mathrm{RSp}}
\newcommand\rspclhit{\Xi_\mathrm{Hit, cl}^\mathrm{RSp}}
\def\thp#1{#1^\mathrm{WL}}
\newcommand\bht{(\partial\calH^2)^{(2)}}
\newcommandx{\ioo}[2]{I_{(#1,#2)}}
\newcommandx{\ioc}[2]{I_{(#1,#2]}}
\newcommandx{\ico}[2]{I_{[#1,#2)}}
\newcommandx{\icc}[2]{I_{[#1,#2]}}
\newcommand{\hs}{{H_\Sigma}}
\newcommand\hsf{H_\Sigma^{[4]}}
\newcommand{\hg}{{H_\Gamma}}
\newcommand\hgf{H_\Gamma^{[4]}}

\makeindex

\begin{document}

\title[The real spectrum compactification of character varieties]{The real spectrum compactification \\ of character varieties:  \\ characterizations and applications\\\vskip1cm
La compactification des vari\'et\'es de caract\`eres \\ par le spectre r\'eel:\\ caract\'erisations et applications}

\author[M.~Burger]{Marc Burger}
\address{Departement Mathematik, ETHZ, R\"amistrasse 101, CH-8092 Z\"urich, Switzerland}
\email{burger@math.ethz.ch}
\author[A.~Iozzi]{Alessandra Iozzi}
\address{Departement Mathematik, ETHZ, R\"amistrasse 101, CH-8092 Z\"urich, Switzerland}
\email{iozzi@math.ethz.ch}
\author[A.~Parreau]{Anne Parreau}
\address{Univ. Grenoble Alpes, CNRS, Institut Fourier, F-38000
  Grenoble, France}
\email{Anne.Parreau@univ-grenoble-alpes.fr}
\author[M.B.~Pozzetti]{Maria Beatrice Pozzetti}
\address{Mathematical Institute, Heidelberg University, Im Neuenheimer feld 205, 69120 Heidelberg, Germany }
\email{pozzetti@mathi.uni-heidelberg.de}
\thanks{Beatrice Pozzetti was partially supported by  DFG project 338644254 and the DFG Emmy Noether grant 427903332.
Marc Burger and Alessandra Iozzi were partially supported by the SNF grant 2-77196-16.
Alessandra Iozzi acknowledges moreover support from U.S. National Science Foundation grants DMS 1107452, 1107263, 1107367 
"RNMS: Geometric Structures and Representation Varieties" (the GEAR Network). 
Marc Burger, Alessandra Iozzi and Beatrice Pozzetti  would like to thank the  National Science Foundation under Grant No. 1440140 
that supported their residence at the Mathematical Sciences Research Institute in Berkeley, California, 
during the semester of Fall 2019 where work on this paper was undertaken.}
\date{\today}

\begin{abstract}  
We announce results on a compactification of general character varieties that has good topological properties 
and give various interpretations of its ideal points. We relate this
to the Weyl chamber length  compactification 
and apply our results to the theory of maximal and Hitchin representations.
\vskip1cm
\noindent
\scshape R\'esum\'e. \upshape 
Cette annonce est un survol de nos r\'esultats concernant la compactification de vari\'et\'es de caract\`eres par le spectre r\'eel.  
Nous relions cette compactification \`a celle obtenue par les fonctions longeurs \`a valeurs dans une chambre de Weyl
et donnons des applications aux repr\'esentations maximales et de Hitchin.
\end{abstract}

\maketitle

\setcounter{tocdepth}{2} 
\tableofcontents   

\section{Introduction}\label{sec:intro}
Let $\bG\leq\SL_n$ be a connected semisimple algebraic group defined over $\QQ$ 
and let $G:=\bG(\RR)$ (or any subgroup of finite index in $\bG(\RR)$).
The character variety $\Xi(\Gamma,G)$ of a finitely generated group $\Gamma$
into $G$ is the quotient by $G$-conjugation of the topological space $\Homred(\Gamma,G)$
of reductive representations, where the latter is equipped with the topology of pointwise convergence.
Our aim here is to give an overview of our results in \cite{BIPP_rsp, BIPP_pos, BIPP_max}: 
we study a compactification of $\Xi(\Gamma,G)$ 
with good topological properties
giving various interpretations of its ideal points (\S~\ref{sec:rsp})
and we apply our study to the theory of maximal  (\S~\ref{sec:max}) and Hitchin representations (\S~\ref{sec:Hpos}).

Ideally a compactification of $\Xi(\Gamma,G)$ should reflect its homological properties 
and its ideal points should reflect ways in which invariants of metric nature associated to reductive representations of $\Gamma$ into $G$
behave as one moves to infinity.
This program can be achieved by using the {\em real spectrum compactification} of $\Xi(\Gamma,G)$ 
and by relating it to the Weyl chamber length compactification obtained by considering Weyl chamber valued length functions.

If $\Gamma=\pis$ is the fundamental group of a connected oriented surface $\Sigma$ with negative Euler characteristic, possibly with non-empty boundary,
we study the action of the mapping class group $\mcgp$  of $\Sigma$ on the real spectrum boundary 
of the components of $\Sp(2n,\RR)$-valued maximal representations.
As a consequence of the theory outlined above, we obtain for example that:
\be
\item all stabilizers for the $\mcgp$-action on the boundary of the components of $\Sp(2n,\RR)$-valued maximal representations
are virtually Abelian, and
\item if $2n\geq4$, there is a non-empty open set of discontinuity for the $\mcgp$-action.
\ee

If $\Sigma$ has no boundary we obtain analogous results concerning
the $\PSL(n,\RR)$-Hitchin components. Furthermore, points in the set
of discontinuity for the $\mcgp$-action on the real spectrum boundary
of a $\PSL(n,\RR)$-Hitchin component or of a component of maximal representations 
give rise to harmonic conformal maps into appropriate buildings.

\medskip
This introduction is by all means not meant to be exhaustive, but only to titillate the curiosity of the reader
and entice his or her to read this announcement.


\section{The real spectrum compactification}\label{sec:rsp}
Echoing the work of Brumfiel \cite{Brum2}, we use techniques and concepts from real algebraic geometry 
in which semialgebraic sets and semialgebraic maps play a fundamental role.  
A semialgebraic set in $\RR^m$ is a finite union of subsets 
defined by finitely many polynomial equalities and inequalities where, in this article, 
the polynomials have coefficients in the field $\qbarr$ of real algebraic numbers.

Such a semialgebraic set admits a {\em real spectrum compactification} (see \cite{Brum89} and \cite[Chapter~7]{BCR}).
For the sake of the reader we give the definition here, although for the understanding of this announcement
it can be taken as a black box.  Let $\calA$ be a commutative ring with unity. 
The real spectrum $\Rspec\calA$ is the set of pairs $(\fp, <)$ 
consisting of a prime ideal $\fp$ in $\calA$ and an ordering $<$ on the field of fractions of $\calA/\fp$. 
A subbasis of open sets for the topology of $\Rspec\calA$ is obtained by declaring open, for each $a\in\calA$, 
the set of pairs $(\fp, <)$ for which $a\mod \fp$ is positive.
Endowed with this topology $\Rspec\calA$ is compact (but not Hausdorff) and its subset of closed points is compact and Hausdorff. 

If $\calA=\qbarr[V]$ is the coordinate ring of an algebraic subset $V\subset \RR^n$ defined over $\qbarr$, 
we let $\rsp V$ be the real spectrum of $\calA$; evaluation on points of $V$ gives an injection $V\hookrightarrow \rsp V$ with open dense image.  
If $W\subset V$ is semialgebraic, we let $\wt W\subset \rsp V$ be the corresponding constructible set, 
namely the subset of $\rsp V$ defined by the same polynomial equalities and inequalities as $W$. 
Then a fundamental theorem \cite[Theorem 7.2.3]{BCR} says 
that the map $W\to \wt W$ is an isomorphism from the Boolean algebra of semialgebraic sets to the one of constructible sets, 
preserving the properties of being open and closed. 
As a consequence if $W\subset V$ is a closed semialgebraic subset and $\rsp W$ is its closure in $\rsp V$ 
then $\rsp W=\wt W$, $W$ is open and dense in $\rsp W$ and the latter is closed and compact.

It is well known (see \cite{Boehm_Lafuente} for an elementary treatment or \cite{RS}) 
that the character variety $\Xi(\Gamma,G)$ is homeomorphic to a closed semialgebraic subset of $\RR^m$ for some $m$.
As such, it admits a real spectrum compactification $\rsp\Xi(\Gamma,G)$ 
whose subset $\rspcl\Xi(\Gamma,G)$ of closed points is the main object of our study.  
Notice that a semialgebraic model of $\rsp\Xi(\Gamma,G)$ depends on
certain choices, but any two sets of choices lead to models that are
related by a canonical semialgebraic isomorphism, and thus to a canonical real spectrum compactification
on which the action of $\Out(\Gamma)$ extends continuously
(see \cite{Brum89} and \cite[Chapter~7]{BCR}).

Let us list at the outset some topological properties of the real spectrum compactification applied to our situation:
\be
\item\label{it:intro1} $\rspcl{\Xi}(\Gamma,G)$ is a compact metrizable $\Out(\Gamma)$-space;
\item\label{it:intro2} not only do the injections
\bqn
\xymatrix{
&\rspcl{\Xi}(\Gamma,G)\ar[dd]^\subset\\
\Xi(\Gamma,G)\quad\ar@{^{(}->}[ur]\ar@{^{(}->}[dr]&\\
&\rsp{\Xi}(\Gamma,G)
}
\eqn
have open dense image
but also they induce bijections at the level of connected components;
\item\label{it:intro4} if $\calV\subset\Xi(\Gamma,G)$ is a contractible connected component that is invariant under an element $\varphi\in\Out(\Gamma)$,
then $\varphi$ has a fixed point in the real spectrum compactification of $\calV$.  
This holds in greater generality if $\varphi$ has a non-vanishing Lefschetz number on $\calV$,
\cite{Brum88B, Brum92}.
\ee

\subsection{A characterization of the boundary of the real spectrum compactification}\label{subsec:RSp}

Our first result gives an interpretation of the points in the boundary
\bqn
\rsp{\partial}\Xi(\Gamma,G):=\rsp{\Xi}(\Gamma,G)\smallsetminus\Xi(\Gamma,G)
\eqn
in terms of representations of $\Gamma$.
To this end we introduce a few objects and concepts.

Any ordered field has a real closure, which is nothing but the the largest algebraic extension to which the order extends, 
and which is unique up to order preserving isomorphism.
Real closed fields can be characterized by saying that every positive element is a square and any odd degree polynomial has a root.
Typical examples are of course $\qbarr$ and $\RR$; if $\FF$ is any real closed field, then $\qbarr\subseteq\FF$, 
but  $\FF\subseteq\RR$ if and only if $\FF$ is Archimedean\footnote{Recall that an order $>$ on a field $\FF$ is {\em Archimedean} if for any $x,y\in\FF$, with $x>y>0$, there exists $n\in \NN$ such that $ny>x$, in which case $\FF$ is an Archimedean ordered field.}.

\begin{ex}\label{ex:orders}  Let $\qbarr(X)$ be the field of rational fractions.  
All Archimedean orders on $\qbarr(X)$ come from injections of $\qbarr(X)$ in $\RR$ via evaluation, 
namely are of the following form:
\be
\item[$b$:]  $f>_b0$ if $f(b)>0$ for $b\in\RR\smallsetminus\qbarr$.
\ee
On the other hand, the following classify all  non-Archimedean orders on $\qbarr(X)$:
\be
\item[$-\infty$:]  $f>_{-\infty}0$ if there is $T>0$ such that $f(t)>0$ for $t\in(-\infty,T)$;
\item[$a_-$:]  $f>_{a_-}0$ if there is $\epsilon>0$ such that $f(t)>0$ for $t\in(a-\epsilon,a)$, for $a\in\qbarr$;
\item[$a_+$:]  $f>_{a_+}0$ if there is $\epsilon>0$ such that $f(t)>0$ for $t\in(a,a+\epsilon)$, for $a\in\qbarr$;
\item[$\infty$:]  $f>_\infty0$ if there is $T>0$ such that $f(t)>0$ for $t\in(T,\infty)$.
\ee
Observe that none of these last orders is Archimedean, 
as for instance 
\bqn
\frac{1}{X-a}>_{a_+} r\qquad\text{ for all }r\in\qbarr
\eqn
and $\frac{1}{X-a}$ is thus an ``infinitely large" element.

We denote by $\ov{\QQ(X)}^\alpha$ for $\alpha\in\{b,-\infty, a_-,a_+,\infty\}$ the corresponding real closures of $\QQ$.

\end{ex}

For any semialgebraic set $W\subset\RR^m$ the ``extension'' $W_\FF\subset\FF^m$ of $W$ to $\FF$ is well defined 
assuming that $\FF$ is any real closed field, \cite[Definition~5.1.2]{BCR}. 
This is defined as the semialgebraic subset of $\FF^m$ defined by the same equalities, inequalities and operations defining $W$.  
In particular, since the group $G$ is a semialgebraic set, we can talk about $G_\FF$ for any real closed field $\FF$.

Let $\FF$ be a real closed field and $\rho\colon\Gamma\to G_\FF$ a homomorphism.  
We say that $\FF$ is {\em $\rho$-minimal} if $\rho$ is not $G_\FF$-conjugate to a representation
with values in $G_\LL$, where $\LL\subset\FF$ is a proper real closed subfield.  
We establish in \cite{BIPP_rsp} that if $\rho:\Gamma\to G_\FF$ is reductive,
then there exists a unique real closed subfield $\FF_\rho\subset\FF$ such that 
$\rho$ is $G_\FF$-conjugate to a representation $\rho'\colon\Gamma\to G_{\FF_{\rho'}}$ 
and $\FF_{\rho'}$ is $\rho'$-minimal.  Finally, 
if $\FF_1$ and $\FF_2$ are two real closed fields,
we say that the representations $(\rho_1,\FF_1)$ and $(\rho_2,\FF_2)$ are {\em equivalent}
if there is an order preserving isomorphism $i\colon\FF_1\to\FF_2$
such that $\rho_2$ and $i\circ\rho_2$ are $G_{\FF_2}$-conjugate.

\begin{thm}[\cite{BIPP_rsp}]\label{thm:bdry}
 Points in $\rsp\partial\Xi(\Gamma,G)$ are in bijective correspondence
with equivalence classes of pairs $(\rho,\FF_\rho)$, where $\rho\colon\Gamma\to G_{\FF_\rho}$ is reductive,
and $\FF_\rho$ is real closed $\rho$-minimal and non-Archimedean.  
Moreover $\FF_\rho$ is of finite transcendence degree over $\qbarr$,
in particular $\FF_\rho$ is countable.
\end{thm}

The property that $(\rho,\FF_\rho)$ represents a point in the boundary of $\rsp\Xi(\Gamma,G)$
is reflected in the fact that $\FF_\rho\supset\qbarr$ is non-Archimedean: 
it contains elements that are infinitely large relative to $\qbarr$.  

\medskip
In addition to the above, the closed points in $\rsp\partial\Xi(\Gamma,G)$ can be further characterized 
in terms of actions on appropriate buildings.
To this end, we will need an explicit construction, inspired by \cite{KTa, Brum2},
of a building associated to $G_\FF$, where $\FF$ is a  non-Archimedean real closed field 
endowed with a non-trivial order-compatible valuation $\nu\colon\FF\to\RR\cup\{\infty\}$.

\begin{remark}
While not all real closed fields admit such a valuation, the existence in the case of a minimal field as in Theorem~\ref{thm:bdry}
is guaranteed by its finite transcendence degree over $\qbarr$ (see \cite[\S~5]{Brum2} or \cite[Chapter~VI, \S~10.3, Theorem~1]{Bourbaki}).
If such an order compatible valuation exists, it is unique up to multiplication by a positive real.
\end{remark}

Without loss of generality, we can assume that the group $\bG$ is invariant under transposition.  
For an example of the concepts to follow, see \S~\ref{subsec:modelBuilding}.
If we regard the set
\bqn
\calP^1(n):=\{A\in M_{n\times n}(\RR):\,A\text{ is symmetric},\, \det(A)=1,\,A\text{ is positive definite}\}
\eqn 
of positive definite matrices of determinant one as the symmetric space associated to $\SL(n,\RR)$,
the $G$-orbit $\calX:=G\cdot\Id$ of $\Id\in\calP^1(n)$ is the symmetric space associated to $G$, 
which is thus a semialgebraic closed subset of $\calP^1(n)$.
Let $\bA<\bG$ be a maximal $\RR$-split torus defined over $\qbarr$, 
which we may assume consisting of symmetric matrices,
let $A:=\bA(\RR)^\circ$ and $\abp\subset A$ a closed (multiplicative) Weyl chamber.  
Then we have the Cartan decomposition $G=K\abp K$, where $K=G\cap \mathrm{O}(n)$, 
by which we define the $\abp$-distance function 
\bqn
\delta\colon \calX\times\calX\to\abp\,.
\eqn
as the unique element in $\delta(x,y)\in\abp$, for $x,y\in\calX$, such that 
\bq\label{eq:A-distance}
g x=\Id\text\quad{ and }\quad g y=\delta(x,y)\Id\
\eq
for some (not unique) $g\in G$.

We say that 
\bqn
N\colon A\to\RR_{\geq1}
\eqn 
is a {\em multiplicative} norm 
if $N$ is 
\be
\item[(N1)]\label{it:norm1} continuous and semialgebraic, invariant under the action of the Weyl group and
\item[(N2)]\label{it:norm2} satisfies the properties that 
\be
\item $N(gh)\leq N(g)N(h)$ for all $g,h\in A$ and
\item $N(g)=1$  if and only if $g=e$.
\ee
\ee
All objects involved being semialgebraic, 
we obtain for every real closed field an $\abp_\FF$-valued distance function 
\bq\label{eq:deltaF}
\delta_\FF\colon \calX_\FF\times\calX_\FF\to\abp_\FF
\eq
defined by the same formula as in \eqref{eq:A-distance}
and an extension $N_\FF\colon A_\FF\to\FF_{\geq1}$ of $N$, which
satisfies properties (N1) and (N2) as well,  \cite[Chapter~5]{BCR}.

Let now $\nu\colon \FF\to\RR\cup\{\infty\}$ be a non-trivial valuation with value group $\Lambda:=v(\FF^\times)$.
We define $\dnf\colon\calX_\FF\times\calX_\FF\to\RR_{\geq0}$ by 
\bq\label{eq:dnf}
\dnf(x,y):=-\nu(N_\FF(\delta_\FF(x,y)))\,.
\eq

\begin{ex}  If $\FF=\RR$ we can take $\nu=-\ln$.  
\end{ex}

If the field $\FF$ is Archimedean, for example in the case of the reals, then $\dnf$ is a a distance on $\calX_\FF$, \cite{Planche}.
This is however never the case if $\FF$ is non-Archimedean, which is the case of interest to us for
representations in the boundary of the real spectrum compactification.  
In this case, $\dnf$ is a pseudodistance and echoing \cite{KTa} 
(see also \cite{Brum2} for $G=\SL(2,\RR)$ and \cite{Bouzoubaa} for $G=\SO(n,1)$), 
we define $\bgf$ as the quotient $\bgf:=\calX_\FF/\sim$, 
where $x\sim y$ if and only if $\dnf(x,y)=0$.
Then $\bgf$ is a metric space with quotient distance still denoted by 
\bqn
\dnf\colon\bgf\times\bgf\to\Lambda\,.
\eqn
It is readily verified that different multiplicative norms define equivalent distances.
As a result the metric completion $\ov\bgf$ is well defined.
We will use later  in our application to the existence of harmonic maps\footnote{Although we will not use it in this paper, 
 notice that $\bgf$ can be endowed with a structure of $\Lambda$-building in the sense of Bennett (see \cite{KTa} for an announcement and
 \cite{Appenzeller} for details).} that $\ov\bgf$ can be endowed also with a $G_\FF$-invariant CAT(0) distance.

The upshot of the above discussion is that if $(\rho,\FF)$ represents a point in $\rsp\partial\Xi(\Gamma,G)$,
there is a canonically associated building $\bgf$, 
and $\Gamma$ acts isometrically both on $\bgf$ and on its completion $\ov\bgf$.

\begin{thm}[\cite{BIPP_rsp}]\label{thm:bdry_cl}  Let $S=S^{-1}$ be a finite generating set of $\Gamma$
and let $E:=S^{2^n-1}\subset\Gamma$.
The following assertions are equivalent:
\be
\item\label{it:bdry_cl1} The representation $(\rho,\FF)$ represents a closed point in $\rsp\partial\Xi(\Gamma,G)$;
\item\label{it:bdry_cl2} The $\Gamma$-action on $\bgf$ does not have a fixed point;
\item\label{it:bdry_cl3} There exists $\eta\in E$ such that $\rho(\eta)$ has positive translation length on $\ov\bgf$.
\ee
\end{thm}
We discuss briefly here the concept of translation length on $G_\FF$, 
which relies on the Jordan projection, both of which will be used in \S~\ref{subsec:RSp->ThP}.

Let $\fa$ be the Lie algebra of $A$ and let ${\fabp}:=\ln(\abp)$
be the model Weyl chamber of $G$.
Then every $g\in G$ is in a unique way a product $g=e_gh_gu_g$, 
where $e_g,h_g$ and $u_g$ pairwise commute, $e_g$ is conjugate into the maximal compact $K<G$,
$u_g$ is unipotent and $h_g$ is the conjugate of a unique element $\exp(X_g)$ with $X_g\in{\fabp}$.
The Jordan projection $J\colon G\to{\fabp}$ is then defined as $J(g):=X_g$.
Given a real closed field $\FF$ endowed with an order preserving valuation $v\colon\FF\to\RR\cup\{\infty\}$,
one can define a Jordan projection $J_\FF\colon G_\FF\to{\fabp}$ taking also values in ${\fabp}$.
Then the translation length
\bqn
\ell_{N_\FF}\colon G_\FF\to\RR_{\geq0}
\eqn
is defined as usual as
\bq\label{eq:tr_l_def}
\ell_{N_\FF}(g):=\inf \left\{\dnf(gx,x):\,x\in\obgf\right\}
\eq
and one can prove that 
\bq\label{eq:tr_l}
\ell_{N_\FF}(g)=\|J_\FF(g)\|_N\,,
\eq
where $\|\,\cdot\,\|_N\colon\fa\to\RR_{\geq0}$ is the norm on the $\RR$-vector space $\fa$ defined
by 
\bqn
\|X\|_N:=\ln (N(\exp(X)))\,.
\eqn

\subsection{The real spectrum and the Weyl chamber length compactification}\label{subsec:RSp->ThP}
We now turn to the relation between the real spectrum compactification, or rather the subset of its closed points, 
and the Weyl chamber length compactification. 

Given a representation $\rho\colon\Gamma\to G_\FF$,
one obtains by composition with the Jordan projection $J_\FF\circ\rho$
a
Weyl chamber valued ``length function'' on $\Gamma$. 
If $\FF=\RR$ such length functions are used by the third author
to define a compactification
$\thp\Xi(\Gamma,G)$ of $\Xi(\Gamma,G)$
such that the character variety $\Xi(\Gamma,G)$ injects into
$\thp\Xi(\Gamma,G)$ as an open dense set \cite{Parreau12}.
Using Theorem~\ref{thm:bdry_cl}~\eqref{it:bdry_cl3} we see from \eqref{eq:tr_l} that if $(\rho,\FF)$ represents a closed point in $\rsp\partial\Xi(\Gamma,G)$,
the ${\fabp}$-valued function $\gamma\mapsto J_\FF(\rho(\gamma))$ is not identically zero and hence defines a point in $\PP({\fabp}^\Gamma)$.

\begin{thm}[\cite{BIPP_rsp}]\label{thm:RSp->ThP}  The above defined map 
\bqn
\rspcl\Xi(\Gamma,G)\to\thp\Xi(\Gamma,G)
\eqn
is continuous, $\Out(\Gamma)$-equivariant and surjective.
\end{thm}

\begin{remark}
If instead of the Jordan projection we take $g\mapsto \ln((\tr(g))^2+2)$ we obtain a compactification which, in the case $G=\SO(n,1)$, is the
Morgan--Shalen compactification of $\Xi(\Gamma, G)$ in $\PP(\RR_{\geq0}^\Gamma)$. The exact analogue of Theorem \ref{thm:RSp->ThP} holds in this context and this
generalizes \cite[Theorem 17]{Mor92} to all semisimple groups.
\end{remark}

The boundary 
\bqn
\thp\partial\Xi(\Gamma,G):=\thp\Xi(\Gamma,G)\smallsetminus\Xi(\Gamma,G)
\eqn 
consists of projective classes $[L]\in\PP({\fabp}^\Gamma)$ of length functions
such that there exists a sequence $([\rho_k]_{k\geq1})$ that diverges in $\Xi(\Gamma,G)$ 
and such that $[J\circ\rho_k]$ converges to $[L]$.
To study divergent sequences in $\Xi(\Gamma,G)$ 
one associates to a symmetric finite generating set $S$ and a homomorphism $\rho\colon\Gamma\to G$
its displacement function 
\bq\label{eq:D_S}
D_S(\rho)(x):=\max_{\gamma\in S}d(\rho(\gamma)x,x)\,,
\eq
where $x\in\calX$.
This is a convex function with respect to $x$ and if $\rho$ is reductive it achieves its minimum
\bq\label{eq:lambdaS}
\lambda_S(\rho):=\min_{x\in \calX}D_F(\rho)(x)\,.
\eq
Then $\lambda_S\colon\Xi(\Gamma,G)\to[0,\infty)$ is a proper function \cite{Parreau12}.
Let $((\rho_k,\RR))_{k\geq1}$ be a sequence with $\lim_{k\to\infty}\lambda_S(\rho_k)=\infty$
and normalized in such a way that 
\bq\label{eq:centered}
\lambda_S(\rho_k)=D_F(\rho_k)(\Id)\,.
\eq
We call {\em centered} a representation that satisfies \eqref{eq:centered}.
Let $\omega$ be a non-principal ultrafilter on $\NN$ and let $\rol$ be the Robinson field
associated to $\omega$ and to the sequence $\pmb{\lambda}=(\lambda_S(\rho_k))_{k\geq1}$. 
Then $\rol$ is a non-Archimedean real closed field endowed with an order-compatible valuation 
$\nu\colon\rol\to\RR\cup\{\infty\}$ that is surjective.
By \cite{Parreau12} the sequence $(\rho_k)_{k\geq1}$ gives rise to a homomorphism 
$\rhol\colon \Gamma\to G_\rol$.
In the same paper it is proven that the ultralimit of the projectivized length functions $[L_{\rho_k}]$ is 
$[L_\rhol]$, and Theorem~\ref{thm:RSp->ThP} suggests that the ultralimit $\rhol$ should define a point in  $\rspcl\Xi(\Gamma,G)$.
Notice however that the Robinson field $\rol$ has uncountable transcendence degree over $\RR$
and hence, by Theorem~\ref{thm:bdry}, it cannot be the minimal $\rhol$-field 
and the representation $(\rhol,\rol)$ cannot represent a point in $\rsp\partial\Xi(\Gamma,G)$.
Nevertheless we can prove that $(\rhol,\rol)$ is in the real spectrum compactification in an appropriate sense and, 
in fact, completely characterize the closed points in $\rspcl\partial\Xi(\Gamma,G)$ in terms of ultralimits, as follows:

\begin{thm}[\cite{BIPP_rsp}]\label{thm:ggt}
Let $\omega$ a non-principal ultrafilter on $\NN$, 
let $((\rho_k,\RR))_{k\geq1}$ be a sequence of centered representations, and let $(\rhol,\rol)$ be its $(\omega,\pmb{\lambda})$-limit.
Then:
\be
\item $\rhol$ is reductive, and
\item if $\FF_{\rhol}$ is the $\rhol$-minimal field, and $(\rhol,\rol)$ is $G_\rol$-conjugate to $(\pi,\FF_{\rhol})$, 
then $(\pi,\FF_{\rhol})$ represents a closed point in $\rsp\partial\Xi(\Gamma,G)$.
\ee
Conversely, any $(\rho,\FF)$ representing a point in $\rspcl\partial\Xi(\Gamma,G)$ arises in this way.
In other words there exists a sequence of centered homomorphisms $((\rho_k,\RR))_{k\geq1}$ as above, 
a non-principal ultrafilter $\omega$ and an order preserving field injection
\bq\label{eq:inj}
i\colon \FF\hookrightarrow\rol
\eq
such that $i\circ\rho$ and $\rhol$ are $G_\rol$-conjugate.
\end{thm}

\begin{remark}\label{rem:afterBIPP_rsp}
\be
\item\label{it:seqInCalV}
In fact the above theorem says more.
If $(\rho_k,\RR)_{k\geq1}$ is a sequence of representations in a closed semialgebraic subset $\calV$ of $\Xi(\Gamma,G)$ 
with closure $\rsp\calV$ in $\rsp\Xi(\Gamma,G)$,
then the $(\omega,\pmb{\lambda})$-limit of $(\rho_k)$ is in $\rsp\partial \calV:=\rsp \calV\smallsetminus \calV$.
Conversely in $(\rho,\FF)$ represents a point in $\rsp\partial\calV$,
then the sequence $(\rho_k)_{k\geq1}$ in Theorem~\ref{thm:ggt} can be taken in $\calV$.

Furthermore
\bqn
J_\FF(\rho(\gamma))=J_{\rol}(\rhol(\gamma))
\eqn
for all $\gamma\in\Gamma$.
\item\label{it:ConverseForNotClosed}
The converse part of Theorem~\ref{thm:ggt} holds also for points
in $\rsp\partial\Xi(\Gamma,G)$ that are not necessarily closed, with
the difference that the Robinson field
$\romu$ is defined by a sequence $\pmb{\mu}$ that grows
faster than $\pmb{\lambda}=(\lambda_S(\rho_k))_{k\geq1}$. 
This accounts for the global fixed point of $\rhomu$ on $\bgromu$ (compare with Theorem~\ref{thm:bdry_cl}).
\ee
\end{remark}
The hypothesis that the sequence $((\rho_k,\RR))_{k\geq1}$ is centered is essential.  
There are indeed examples of sequences of non-centered reductive representations whose ultralimit is not reductive.

Notice that if $\FF_j$, $j=1,2$, are real closed fields with order-compatible valuations $\nu_j$,
then any field injection $i\colon\FF_1\to\FF_2$ compatible with $\nu_1$ and $\nu_2$ induces an embedding
\bqn
\calB_{\FF_1}\hookrightarrow\calB_{\FF_2}
\eqn
that is isometric with respect to the distances associated to a multiplicative norm. 
The building $\bgrol$ is an affine $\RR$-building
with a complete CAT(0) metric for which it is isometric to the asymptotic cone of the sequence of metric spaces 
$(\calX,\Id,\frac{d}{\lambda_S(\rho_k)})_{k\geq1}$ \cite{Parreau12}. 

\begin{prop}[\cite{BIPP_rsp}]\label{prop:isom_emb}
The embedding  
\bqn
\obgf\hookrightarrow \bgrol
\eqn
coming from the field injection $i$ in \eqref{eq:inj} identifies $\obgf$ with a subspace of $\bgrol$ 
that is totally geodesic and complete with respect to the CAT(0) metric on $\bgrol$.
\end{prop}

This theorem allows us to use the geometry of symmetric spaces to study the actions on buildings
associated to points in the real spectrum boundary.
It is particularly useful when applied to a connected component $\calV\subset\Xi(\Gamma,G)$
consisting of representations of geometric interest.

\subsection{Harmonic maps into buildings}\label{subsec:harm}
The above results can be used to show the existence of equivariant harmonic maps.

Let $(Y,d)$ be a general CAT(0) space, $\rho:\Gamma\to\Iso(Y)$ a homomorphism and $S$ a symmetric finite generating set.
The isometric $\Gamma$-action on $Y$ is {\em KS-proper} 
if for every $T\geq0$ the set 
\bqn
\{y\in Y:\,D_S(\rho)(y)\leq T\}
\eqn 
is bounded, \cite{KSch}, where $D_S$ is as in \eqref{eq:D_S}.
This property is independent of the choice of the finite generating set $S\subset \Gamma$.
Properness implies the -- usually -- weaker property that the action is non-parabolic,
that is $\Gamma$ does not fix a point in the visual boundary of $Y$.
If closed balls are compact, the properties are equivalent.  

In our context so far, we say that a representation $\rho:\Gamma\to G$ is {\em non-parabolic}
if the corresponding $\Gamma$-action on the associated symmetric space $\calX$ is non-parabolic.

Using the full strength of Theorem~\ref{thm:ggt} and knowing from Proposition~\ref{prop:isom_emb}
that $\obgf$ is a CAT(0) space, we can prove the following:

\begin{thm}[\cite{BIPP_rsp}]\label{thm:proper}  Let $\calV\subset\Xi(\Gamma,G)$ be a closed semialgebraic set
consisting of non-parabolic representations.  Then for any $(\rho,\FF)$ in $\rsp\partial\calV$
the corresponding action on $\obgf$ is KS-proper.
\end{thm}

Then applying \cite[Theorem~2.1.3]{KSch} we obtain:

\begin{cor}[\cite{BIPP_rsp}]\label{cor:harm}  Let $\Gamma=\pi_1(M)$, where $M$ is a compact Riemannian manifold
and let $\calV\subset\Xi(\Gamma,G)$ be a closed semialgebraic set consisting of non-parabolic representations.
Then for any $(\rho,\FF)\in\rsp\partial\calV$ there exists a $\pi_1(M)$-equivariant Lipschitz harmonic map $\wt M\to\obgf$.
\end{cor}

Notice that the above result applies also to actions with a global fixed point, in which case the equivariant harmonic map
is obviously constant.

\subsection{Accessibility}\label{subsec:access}
We turn to another application of the link between the real spectrum and the Weyl chamber length compactifications 
in Theorem~\ref{thm:RSp->ThP}.  Recall that already in the case in which $\Sigma$ is a compact surface of genus $g\geq2$,
$\Gamma=\pi_1(\Sigma)$ and $G=\PSL(2,\RR)$, the topological space $\thp\Xi(\Gamma,G)$ is rather badly behaved 
as the closure of the connected components of Euler number $1\leq e\leq 2g-3$ contain the Thurston
boundary as a nowhere dense subset,  \cite{Wolff}.  
However using the real spectrum compactification and its general algebraic properties allows us
to derive some nice topological properties of the Weyl chamber length compactification.
\begin{thm}[\cite{BIPP_rsp}]\label{thm:semialg}  Let $V\subset\RR^n$ be a semialgebraic set, 
let $\Phi\colon V\to[0,\infty)$ be a proper semialgebraic continuous function
and let $\alpha\in\rspcl\partial V$ be a closed point in the real spectrum boundary of $V$.
Then there exists a continuous and locally semialgebraic section of $\Phi$ with endpoint $\alpha$, namely there exists $T>0$ and 
\bqn
\sigma_\alpha\colon[T,\infty)\to V
\eqn
continuous and locally semialgebraic, such that
\be
\item  $\Phi\circ\sigma_\alpha(t)=t$;
\item $lim_{t\to\infty} \sigma_\alpha(t)=\alpha$.
\ee
\end{thm} 
  
As a consequence we have that, although badly behaved, 
every connected component of the Weyl chamber length compactification is pathwise connected.

\begin{thm}[\cite{BIPP_rsp}]\label{thm:path1}  
For every $[L]\in\thp\partial\Xi(\Gamma,G)$ there exists a continuous and locally semialgebraic path
\bqn
\ba{}
[1,\infty)&\to\Homred(\Gamma,G)\\
t\quad&\longmapsto\quad\quad\rho^t
\ea
\eqn
such that 
\bqn
\lim_{t\to\infty}\frac{J(\rho^t(\gamma))}{t}=L(\gamma)
\eqn
for every $\gamma\in\Gamma$.

In particular 
\bqn
\lim_{t\to\infty}[\rho^t]=[L]
\eqn
and every connected component of the Weyl chamber length compactification is path connected.
\end{thm}

Recall that we introduced in \eqref{eq:lambdaS} the function $\lambda_S(\rho)$, 
that is the minimal joint displacement of the generating set $S$.
Another way to measure the displacement of a representation $\rho$ is,
for any finite set $E\subset\Gamma$, to consider the maximal translation length of an element in $E$
\bqn
\ell_E(\rho):=\max_{\gamma\in E}\ell(\rho(\gamma))\,.
\eqn
There is the obvious inequality $\ell_E(\rho)\leq\lambda_E(\rho)$.
Taking as input the work in \cite{RS}, we prove:

\begin{thm}[\cite{BIPP_rsp}]\label{thm:path2}  
Let $S=S^{-1}$ be a finite generating set and let $E=S^{2^n-1}\subset\Gamma$.
Then on $\Xi(\Gamma,G)$ the functions $\rho\mapsto\lambda_S(\rho)$ and $\rho\mapsto\ell_E(\rho)$ are Lipschitz equivalent.
\end{thm}

\begin{cor}[\cite{BIPP_rsp}]\label{cor:path}  
Let $[L]\in\thp\partial\Xi(\Gamma,G)$ and let $(\rho^t)_{t\geq1}$ be as in Theorem~\ref{thm:path1}.
Then the functions $\lambda_S(\rho^t)$ and $\ell_E(\rho^t)$ are Lipschitz equivalent to $t$.
\end{cor}

\subsection{Applications to cyclic Higgs bundles}\label{subsec:higgs}
We can give sufficient conditions on a path of representations $(\rho^t)$ 
guaranteeing convergence to a unique closed point in the real spectrum compactification
in terms of asymptotic expansions of the functions $t\mapsto\tr(\rho^t(\gamma))$.
More precisely, assume for the sake of simplicity that $\calV\subset\Xi(\Gamma,G)$
is a closed semialgebraic subset such that the representations in $\calV$ are determined by their character\footnote{Important examples of components $\calV$ with the above properties are the Hitchin component if $G=\SL(n,\RR)$
and the set of maximal representations if $G=\Sp(2n,\RR)$, \cite{BIPP_max}.}.
Then we have:
\begin{thm}[\cite{BIPP_rsp}]\label{thm:asexp}  Let $\calV$ be as above, 
$t\mapsto\rho^t$ be a continuous path of representations with image in  $\calV$
and $S\subset\Gamma$ a symmetric finite generating set.
Assume that for each $\gamma\in S^{2^n-1}$ there exists an asymptotic expansion
\bq\label{eq:2.15}
\tr(\rho^t(\gamma))=\sum_{j=1}^{n_\gamma} a_{\gamma,j}t^{\alpha_{\gamma,j}}+R_\gamma(t)\,,
\eq
where for all $k\in\ZZ$, $t^kR_\gamma(t)$ is bounded as a function of $t\geq1$.
If there exists $\gamma$ with $\alpha_{\gamma,1}>0$, then the path has a unique limit in $\rsp\partial\calV$
as $t\to\infty$.
\end{thm}
The motivation for this general result arose from conversations with
J.~Loftin, A~Tamburelli and M.~Wolf in the context where $\Gamma=\pi_1(\Sigma)$ is
the fundamental group of a closed orientable surface $\Sigma$ of
negative Euler characteristic and $\calV$ is the $\SL(3,\RR)$-Hitchin
component. In fact under the Labourie--Loftin parametrisation of
$\calV$ by the bundle $\calQ_3(\Sigma)$ of cubic differentials over
Teichm\"uller space, the sequence of representations $\rho^t$
corresponding to $s\cdot q$, where $\ln t=s^{1/3}$ and
$q\in\calQ_3(\Sigma)$ should satisfy \eqref{eq:2.15}, 
at least for geodesic paths that avoid the zeroes of q, \cite{Loftin_Tamburelli_Wolf}.
Thus this path of representations has a unique limit in the real spectrum boundary.

\section{Maximal representations}\label{sec:max}
Let $\Sigma$ be a connected compact oriented surface possibly with boundary $\partial\Sigma$,
and assume that $G:=\bG(\RR)^\circ$ is of Hermitian type, 
that is the associated symmetric space is Hermitian.
Then the  set $\Ximax(\Sigma,G)\subset\Xi(\pis,G)$ consisting of $G$-conjugacy classes of maximal representations
is a closed semialgebraic subset, \cite{BIW}, 
and if $\partial\Sigma=\varnothing$, it is a union of connected components of $\Xi(\pis,G)$.
Notice that if $\partial\Sigma\neq\varnothing$, then $\pis$ is a free group, 
but the notion of maximality depends heavily on $\Sigma$ and not only on its fundamental group.

The closure $\rsp\Ximax(\Sigma,G)$ of the
character variety $\Ximax(\Sigma,G)$ of maximal representations  in
$\rsp\Xi(\pis,G)$ is then
 a compactification of $\Ximax(\Sigma,G)$, that is a compact space containing
$\Ximax(\Sigma,G)$ as an open dense subset.
%
We will now present results that, building on the general theory outlined in \S~\ref{sec:intro}, show 
how the geometric and dynamical properties of maximal representations are reflected in the behavior of the objects
associated to points in the real spectrum boundary
\bqn
\rsp\partial\Ximax(\Sigma,G):=\rsp\Ximax(\Sigma,G)\smallsetminus\Ximax(\Sigma,G)\,.
\eqn
Notice that since $\Ximax(\Sigma,G)$ is closed in $\Xi(\pis,G)$, then 
\bqn
\rsp\partial\Ximax(\Gamma,G)\subset\rsp\partial\Xi(\pis,G)\,.
\eqn

While our methods apply in full generality, that is for any  Lie group of Hermitian type,
for the sake of concreteness we will restrict to the case of $G=\Sp(2n,\RR)$.

\subsection{Definition and characterization of maximal representations}\label{subsec:def_char}
Our objective in this section is to establish two characterizations
 of the points in $\rsp\Ximax(\Sigma,\Sp(2n,\RR))$.
The first is of cohomological nature.
To this purpose, recall that according to Theorem~\ref{thm:bdry},
points in $\rsp\Xi(\pis,\Sp(2n,\RR))$ are represented by $\rho\colon\pis\to\Sp(2n,\RR)$, where $\FF$ is $\rho$-minimal.  
We show here that points in $\rsp\Ximax(\Sigma,\Sp(2n,\RR))$ are representations of $\pis$ into $\Sp(2n,\FF)$ 
that are {\em maximal} according to a definition that parallels the one of maximal representations into $\Sp(2n,\RR)$.
Given any ordered field $\KK$, we will equip $\ktn$ with the standard symplectic form
\bqn
\left\langle\bpm x_1\\y_1\epm, \bpm x_2\\y_2\epm\right\rangle:={}^\mathrm{t}x_1y_2-{}^\mathrm{t}x_2y_1\,.
\eqn
Let $\calL(\ktn)$ be the Grassmannian of Lagrangian subspaces of $\ktn$
and let $\tau\colon\calL(\ktn)^3\to\ZZ$ be the Maslov cocycle, which is defined as the signature of the quadratic form
\bqn
\ba
q\colon\ell_1\times\ell_2\times\ell_3&\longrightarrow \qquad\qquad\KK\\
(x_1,x_2,x_3)\,&\mapsto \langle x_1,x_2\rangle+\langle x_2,x_3\rangle+\langle x_3,x_1\rangle
\ea
\eqn
for $\ell_1,\ell_2,\ell_3\in\calL(\ktn)$.
Since it takes values in $[-n,n]$, we obtain in a standard way a bounded cohomology class $\kkb\in\hb^2(\Sp(2n,\KK),\RR)$.  
We use this class to associate to a homomorphism $\rho\colon\pis\to\Sp(2n,\KK)$ its {\em Toledo invariant} $\T(\rho)\in\RR$, \cite[Definition~3.1]{BIW},
which satisfies the inequality
\bqn
|\T(\rho)|\leq2n|\chi(\Sigma)|\,.
\eqn

\begin{defi}\label{defi:max}  A representation $\rho\colon\pis\to\Sp(2n,\KK)$ is {\em maximal}  if 
\bqn
\T(\rho)=2n|\chi(\Sigma)|\,.
\eqn
\end{defi}
The second characterization will use the concept of maximal framing introduced in \cite{BP}.
We say that a triple $(\ell_1,\ell_2,\ell_3)\in\calL(\ktn)^3$ of Lagrangians is {\em maximal} if
$\tau(\ell_1,\ell_2,\ell_3)=n$.  Endow $\mathring{\Sigma}$ with a finite area complete hyperbolic structure and 
identify $\pis$ with a lattice in $\PSL(2,\RR)$ via the holonomy representation,
so that it acts on the hyperbolic plane $\calH^2$ and on its circle at infinity $\partial\calH^2$.

\begin{defi}\label{defi:framing}  A {\em maximal framing} of a homomorphism $\rho\colon\pis\to\Sp(2n,\KK)$ is a pair $(X,\varphi)$
consisting of a non-empty $\pis$-invariant subset $X\subset\partial\calH^2$ and an equivariant map 
$\varphi\colon X \to\calL(\ktn)$ such that the image of any positively oriented triple in $X^3$ is a maximal triple in $\calL(\ktn)^3$.
\end{defi}
Let $\hs\subset\partial\calH^2$ denote the set of fixed points of hyperbolic elements in $\pis$.

\begin{thm}[\cite{BIPP_rsp}]\label{thm:max=fr}  Let $\rho\colon\pis\to\Sp(2n,\FF)$ be a homomorphism,
where $\FF$ is a real closed field.  The following are equivalent:
\be
\item\label{it:max=fr2} The representation $\rho$ is maximal;
\item\label{it:max=fr4} The representation $\rho$ admits a maximal framing defined on $\hs$.
\ee
\end{thm}

\begin{remark}\label{rem:max=f3}In our applications the central property of maximal representations is
the existence of a maximal framing defined on $\hs$.
However there are cases in which there is an interest in considering maximal framing defined 
on appropriate subsets of $\partial\calH^2$.  
For example in \cite{BP} the authors showed that if $\rhol\colon\pis\to\Sp(2n,\rol)$ is the ultralimit of a sequence $(\rho_k)_{k\geq1}$
of maximal representations into $\Sp(2n,\RR)$, then $\rhol$ admits a maximal framing defined on the full boundary $\partial\calH^2$.
At the other end of the spectrum, if $\FF$ is $\rho$-minimal, it is countable and 
hence one cannot expect the maximal framing to be defined on the whole of $\partial\calH^2$.  
  Another reason for the need of flexibility in the definition of a
  maximal framing is that the Fock-Goncharov coordinates on the
  $\Sp(2n,\RR)$-Hitchin component  of a surface $\Sigma$ with
  boundary  lead to maximal framings that are defined on the
subset of $\partial\calH^2$ consisting of the cusps of $\pis$.
We proved that an equivalent condition to those in Theorem~\ref{thm:max=fr} is that 
there exists a $\pis$-invariant subset $X\subset\partial\calH^2$ and an order preserving field injection
$i\colon\FF\to\LL$, where $\LL$ is real closed, such that $i\circ\rho$ admits a maximal framing defined on $X$;
\end{remark}

\begin{ex}\label{ex:maxrep} Let $\Sigma_{0,3}$ be a pair of pants and let
\bqn
\pi_1(\Sigma_{0,3})=\langle c_1,c_2,c_3:\,c_3c_2c_1\rangle
\eqn
be a presentation of its fundamental group.
Define $\rho:\pi_1(\Sigma_{0,3})\to\Sp(4,\QQ(X))$ by
\bqn
\ba
\rho(c_1)
&:=\bpm
1&4X&0&0\\
0&1&0&0\\
2&4X&1&0\\
-4X&2&-4X&1
\epm
\\
\rho(c_2)
&:=\bpm
1&\frac{1}{X}&-2&-\frac{1}{X}\\
0&1&\frac{1}{X}&-2\\
0&0&1&0\\
0&0&-\frac{1}{X}&1
\epm\,.
\ea
\eqn
Composing $\rho$ with the group homomorphism induced by the fields embedding
$\QQ(X)\hookrightarrow\ov{\QQ(X)}^\alpha$, for $\alpha\in\{-\infty, a_-,a_+,\infty\}$ (see Example~\ref{ex:orders})
we obtain representations $\rho_\alpha$ such that $\rho_{0_+}$, $\rho_{a_-}$ and $\rho_{a_+}$ (both for $a\in\qbarr$, $a>0$), and $\rho_\infty$ are maximal.
This rests on Theorem~\ref{thm:max=fr}\,(\ref{it:max=fr4}) and the following facts:
\be
\item  Denoting by $\rho^t:\pi_1(\Sigma_{0,3})\to\Sp(4,\RR)$ the specialization of $\rho$ to $X=t$,
it follows from the explicit parametrization of $\Ximax(\Sigma_{0,3},\Sp(2n,\RR))$ in \cite{Strubel}
that $\rho^t$ is maximal for every $t>0$.
\item Denoting by $\varphi_t\colon \hg\to\mathcal{L}(\RR^4)$ the maximal framing associated to $\rho^t$ by \cite{BIW}, 
the assignment $t\mapsto\varphi_t(\xi)$ is continuous and semialgebraic  for every $\xi\in\hg$.  
The desired maximality for $\rho_{0_+}\colon \pi_1(\Sigma_{0,3})\to\Sp(4,\ov{\QQ(X)}^{0_+})$ 
is a consequence of the maximality of $\rho^t$ for all $t>0$ and the identification of $\ov{\QQ(X)}^{0_+}$
with the field of Puiseux series (with $\qbarr$-coefficients) that are
convergent in some interval $(0,\epsilon)$, for $\epsilon>0$.
\ee
An analogous argument works for the orders $>_\infty, >_{a_-}$ and $>_{a_+}$, for $a>0$.
\end{ex}

\subsection{Description of the real spectrum compactification of $\Ximax(\Sigma,\Sp(2n,\RR))$}
While in general it is difficult to give a concrete description of the minimal field of a point in the real spectrum compactification,
one of the first applications of Theorem~\ref{thm:max=fr} is the
determination of the minimal field $\FF_\rho$ of definition of a
maximal representation $\rho$.
Observe that if $\rho\colon\pis\to\Sp(2n,\FF)$ is any homomorphism, not necessarily maximal, $\FF$ contains the field $\QQ(\tr(\rho))$
generated over $\QQ$ by the set $\{\tr(\rho(\gamma)):\,\gamma\in\pis\}$ of all traces.

\begin{thm}[\cite{BIPP_max}]
  \label{thm:tr}
Let $\rho\colon\pis\to\Sp(2n,\FF)$ be a maximal representation.
Then the minimal  field $\FF_\rho$ of definition of $\rho$ is
the real closure of the field $\QQ(\tr(\rho))$ of traces of $\rho$.
\end{thm}

The next result will ensure that the real spectrum boundary of $\Ximax(\Sigma,\Sp(2n,\RR))$
consists exactly of all maximal representations $(\rho,\FF)\in\rsp\partial\Xi(\pis,\Sp(2n,\RR))$,
where $\FF$ is the real closure of $\QQ(\tr(\rho))$.  
In fact, on the one hand one can prove that any representation in $\rsp\partial\Ximax(\Sigma,\Sp(2n,\RR))$
admits a maximal framing.  
On the other, one can see that if $(\rho,\FF)\in\rsp\partial\Xi(\pis,\Sp(2n,\RR)$,
then $(\rho,\FF)$ is an appropriate limit of a sequence of representations 
that can be taken to be maximal if $\rho$ was.  
(This follows from a version of Theorem~\ref{thm:ggt} for ultraproducts, of which Theorem~\ref{thm:ggt} is a consequence.)


\begin{thm}[\cite{BIPP_max}]\label{thm:bdryMax}
There is a bijective correspondence between points in the real spectrum compactification 
$\rsp\Ximax(\Sigma,\Sp(2n,\RR))$ and equivalence classes $(\rho,\FF)$,
where $\rho\colon\pis\to\Sp(2n,\FF)$ is maximal and $\FF$ is the real closure of the field $\QQ(\tr(\rho))$ of traces of $\rho$.

In addition $(\rho,\FF)$ represents a closed point if and only if $\FF$ is Archimedean over the ring $\QQ[\tr(\rho)]$
generated by the set of traces of $\rho$.
\end{thm}

\begin{ex}\label{ex:rho_closed} (Example \ref{ex:maxrep} continued)
  Since
\bqn
\tr(\rho(c_1^{-1}c_3)^2)=-256X^2+320-\frac{16}{X^2}\,,
\eqn
it follows from Theorem~\ref{thm:bdryMax} that both $(\rho_{0_+},\ov{\QQ(X)}^{0_+})$ and $(\rho_\infty,\ov{\QQ(X)}^\infty)$
are closed points in $\rsp\partial\Ximax(\Sigma_{0,3},\Sp(4,\RR))$.

On the other hand, since the ring of traces of $\rho$ is contained in $\QQ[\frac{1}{X},X]$, neither of the points
$(\rho_{a_-},\ov{\QQ(X)}^{a_-})$ or $(\rho_{a_+},\ov{\QQ(X)}^{a_+})$ is closed for $a>0$.  
One can show that the unique closed point in their closure \cite[Proposition~7.1.24]{BCR} 
is the specialization $\rho^a\colon\pi_1(\Sigma_{0,3})\to\Sp(4,\RR)$ of $\rho$ at $X=a$.
\end{ex}

\subsection{Small stabilizers}\label{subsec:smallstab}
We can give information on the stabilizers of the $\pis$-action on the CAT(0) space $\ov\bsp$
associated to a point $(\rho,\FF)$ in the boundary of the real spectrum compactification of $\Ximax(\Sigma,\Sp(2n,\RR))$ (Proposition~\ref{prop:isom_emb}).
Namely:

\begin{thm}[\cite{BIPP_max}]\label{thm:smallstab}  Let $(\rho,\FF)\in\rsp\partial\Ximax(\Sigma,\Sp(2n,\RR))$
and let $\ov\bsp$ be the associated CAT(0) space. 
Then the stabilizer of a germ of a geodesic segment in $\ov\bsp$ is an elementary subgroup of $\pis$.
\end{thm}

The proof uses the following mechanism of independent interest:

 \begin{remark} (Reduction modulo an order convex subring)
Let $(\rho,\FF)$ be a point in $\rsp\partial\Ximax(\Sigma,\Sp(2n,\RR))$.  
Let $M\subset\FF$ be an order convex subring, 
in particular a valuation ring, and $k$ the residue field (which is real closed as well).
Assume that $\rho(\pis)<\Sp(2n,M)$ and let $\sigma\colon\pis\to\Sp(2n,k)$
be the homomorphism obtained by composing $\rho$ with the reduction map $\Sp(2n,M)\to\Sp(2n,k)$.  
Let $(\sigma,k_\sigma)$ denote the corresponding point in $\rsp\Xi(\pis,\Sp(2n,\RR))$.
It follows directly from the definition of the topology
on the real spectrum that the point  $(\sigma,k_\sigma)$ is in the closure
of the point $(\rho,\FF)$; since $\rsp\Ximax(\Sigma,\Sp(2n,\RR))$ is closed and consists of maximal
representations (Theorem~\ref{thm:bdryMax}), this implies that  
$(\sigma,k_\sigma)$ is maximal as well.
\end{remark}

For the proof of Theorem~\ref{thm:smallstab}, 
we use the converse point of Theorem~\ref{thm:ggt}  (see also Remark~\ref{rem:afterBIPP_rsp})
and, using Proposition~\ref{prop:isom_emb}, reduce ourselves to the case of an ultralimit $\rhol$ of centered maximal representations.
Let $\ol$ be the valuation ring coming from the canonical valuation of the Robinson field
and $\kol$ the residue field.
To every germ $c$ of geodesic segment starting at the basepoint of $\bspkol$ corresponds a proper parabolic subgroup
$P<\Sp(2n,\kol)$, namely
\bqn
\Stab_{\Sp(2n,\kol)}(c)=\{g\in\Sp(2n,\ol):\,\text{the image of }g\text{ in }\Sp(2n,\kol)\text{ is in }P\}\,.
\eqn
If the stabilizer $\Gamma_1$ of $c$ in $\pis$ is non-elementary,
we may find a subsurface $\Sigma_1\subset\Sigma$
with $\pi_1(\Sigma_1)<\Gamma_1$.  
Notice that the characterization in Theorem~\ref{thm:max=fr} of maximal representations in terms of maximal framings
defined on $\hs$ implies that if $\Sigma_1\subset\Sigma$ is a subsurface and $\rho\colon\pis\to\Sp(2n,\FF)$ 
is maximal, its restriction to $\pi_1(\Sigma_1)$ is maximal as well. 
By restricting $\rhol$ to $\pi_1(\Sigma_1)$ and composing with the reduction map
\bqn
\Sp(2n,\ol)\longrightarrow\Sp(2n,\kol)
\eqn
we obtain a maximal representation with image  contained in a parabolic subgroup.
But then the characterization of maximal representations in terms of maximal framings implies that this is a contradiction.

\subsection{A model for the building}
\label{subsec:modelBuilding}
Let $\FF$ be a real closed field and $\Sym_n(\FF)$ the set of symmetric matrices with coefficients in $\FF$.
If $\FF=\RR$, as a model for the symmetric space associated to $\Sp(2n,\RR)$ we take the Siegel upper half space
\bqn
\calS^n:=\{Z=X+iY:\,X,Y\in\Sym_n(\RR), Y>>0\}
\eqn
on which $\Sp(2n,\RR)$ acts as
\bqn
g Z:=(AZ+B)(CZ+D)^{-1}\,,
\eqn
where $g=\bpm A&B\\C&D\epm\in\Sp(2n,\RR)$.
In general $\Sp(2n,\FF)$ acts on 
\bqn
\calS^n_\FF:=\{X+iY:\,X,Y\in\Sym_n(\FF), Y>>0\}
\eqn 
by the same formula and, since $\FF$ is real closed, it shares the same transitivity properties with the $\Sp(2n,\RR)$-action on $\calS^n$.

As a maximal $\RR$-split torus defined over $\QQ$ we take the subgroup of $\Sp(2n,\CC)$ consisting of diagonal matrices.
We have then 
\bqn
A_\FF=\left\{\bpm D&0\\0&D^{-1}\epm:\,D=\diag(a_1,\dots,a_n),\,a_i\in\FF_{>0},\,1\leq i\leq n\right\}
\eqn
and as closed Weyl chamber we choose
\bqn
\abp_\FF=\left\{\bpm D&0\\0&D^{-1}\epm\in A:\,a_1\geq\dots\geq a_n\geq1\right\}\,.
\eqn
We consider then the following multiplicative norm $N_\FF\colon A_\FF\to\FF_{\geq1}$
\bq\label{eq:nf_ex}
N_\FF\left(\bpm D&0\\0&D^{-1}\epm\right):=\prod_{i=1}^n\max\{a_i,a_i^{-1}\}
\eq
and a non-trivial valuation $\nu\colon\FF\to\RR\cup\{\infty\}$.
As in \S~\ref{subsec:RSp} 
we consider the Weyl chamber valued distance function $\delta_\FF\colon\calS^n_\FF\times\calS^n_\FF\to\abp_\FF$ defined in \eqref{eq:A-distance}
and the distance or pseudodistance
\bqn
\dnf(Z_1,Z_2)=-\nu(N_\FF(\delta_\FF(Z_1,Z_2)))\,.
\eqn 
For example if $\FF$ is non-Archimedean, this is not a distance:  
in fact, since the valuation of a real number is always zero, then $\dnf(Z_1,Z_2)=0$ 
whenever $\delta_\FF(Z_1,Z_2)\in\abp\cap \abp_\FF$,
since $N(A)\subset\RR_{\geq1}$.

As in \S~\ref{subsec:RSp}, we let $\bsp$ denote the quotient of $\calS^n_\FF$ by the $\dnf=0$ relation.  
Then, as in \eqref{eq:tr_l_def} and \eqref{eq:tr_l} the translation length of an element $g\in\Sp(2n,\FF)$ is given by
\bqn
\ell_{N_\FF}(g)=\sum_{j=1}^n-\nu(|\lambda_j|)\,,
\eqn
where $\lambda_1,\dots,\lambda_n\lambda_n^{-1},\dots,\lambda_1^{-1}$ are the eigenvalues of $g$ counted with multiplicity
and ordered in such a way that
\bqn
|\lambda_1|\geq\dots\geq|\lambda_n|\geq1\,.
\eqn
Here $|\lambda|\in\FF_{\geq0}$ is the square root of the norm of $\lambda\in\FF[i]$.
If $\FF=\RR$, the same formula holds for the translation length of $g$ with $\nu$ replaced by $-\ln$.

\subsection{The geodesic current associated to a maximal representation}
We endow the surface $\mathring{\Sigma}$  with a complete hyperbolic metric of finite area. 
We recall that a geodesic current on $\mathring{\Sigma}$ is a flip-invariant $\pis$-invariant positive Radon measure on the space 
\bqn
\bht:=\{(x,y)\in(\partial\calH^2)^2:\,x\neq y\}
\eqn 
of distinct pairs of points in $\partial\calH^2$.  Basic examples are the current
\bqn
\delta_c:=\sum_{\eta\in\pis/\langle\gamma\rangle}\delta_{\eta(\gamma_-,\gamma_+)}\,,
\eqn
where $(\gamma_-,\gamma_+)$ is the axis of a hyperbolic element $\gamma$ representing a closed geodesic $c\subset\mathring{\Sigma}$,
and the Liouville current $\calL$ which is the unique (up to scaling) $\PSL(2,\RR)$-invariant Radon measure on $\bht$.
Let $i(\mu,\lambda)$ be the Bonahon intersection of two geodesic currents $\mu,\lambda$ 
and recall that 
\bqn
i(\calL,\delta_c)=\ell(c)\,,
\eqn
where $\ell(c)$ is the hyperbolic length of $c$.

\begin{thm}[\cite{BIPP_pos}]\label{thm:repr_curr}  Let $\rho\colon\pis\to\Sp(2n,\FF)$ be a maximal representation,
where $\FF$ is a real closed field endowed with a valuation 
\bqn
\nu\colon\FF\to\RR\cup\{\infty\}
\eqn
that we assume to be equal to $-\ln$ if $\FF=\RR$.
Then there is a canonical geodesic current $\mu_\rho$ on $\mathring{\Sigma}$ such that 
\bqn
i(\mu_\rho,\delta_c)=\ell_{N_\FF}(\rho(\gamma))
\eqn
for every closed geodesic $c\subset\mathring{\Sigma}$ and hyperbolic element $\gamma$ representing $c$.
\end{thm}


Theorem~\ref{thm:repr_curr} has been shown by \cite{Martone_Zhang} if $\Sigma$ has no boundary and $\FF=\RR$. 
Their proof uses crucially the fact that in this case the maximal framing is defined and continuous on $\partial\calH^2$, 
and a classical construction of Hamendst\"adt that associates a geodesic current to a continuous crossratio, \cite{H5}.  
In our situation crossratios will only be defined on $\hs$ and we cannot expect any continuity property.  
In fact it is only positivity that plays a crucial role and allows us to place ourselves in an abstract framework that we now describe.  

Let $\Gamma<\PSL(2,\RR)$ be a torsion free lattice and let $\hg\subset\partial\calH^2$ 
be the set of fixed points of hyperbolic elements in $\Gamma$.
A {\em positive crossratio} is a $\Gamma$-invariant function 
\bqn{}
[\,\cdot\,,\,\cdot\,,\,\cdot\,,\,\cdot\,]\colon\hgf\to[0,\infty)
\eqn
that is defined on the set $\hgf$ of positively oriented quadruples in $\hg$ and satisfies the following properties:
\be
\item\label{it:cr1}  $[x_1,x_2,x_3,x_4]=[x_3,x_4,x_1,x_2]$, and
\item\label{it:cr2} $[x_1,x_2,x_4,x_5]=[x_1,x_2,x_3,x_5]+[x_1,x_3,x_4,x_5]$
\ee
whenever defined (see Remark~\ref{rem:mu}).
Then if $\gamma\in\Gamma$ is a hyperbolic element with fixed points $\{\gamma_-,\gamma_+\}\subset \hg$,
the period $\per(\gamma)$ of $\gamma$ with respect to $[\,\cdot\,,\,\cdot\,,\,\cdot\,,\,\cdot\,]$ is defined by
\bqn
\per(\gamma):=[\gamma_-,x,\gamma x,\gamma_+]\,,
\eqn
where $x\in \hg$ is any point such that  $(\gamma_-,x,\gamma x,\gamma_+)\in\hgf$.
Recall from classical hyperbolic geometry that for the standard crossratio 
the period of a hyperbolic element is its translation length.

\begin{thm}[\cite{BIPP_pos}]\label{thm:per}  Given a positive crossratio $[\,\cdot\,,\,\cdot\,,\,\cdot\,,\,\cdot\,]$ on $\hgf$,
there exists a geodesic current $\mu$ on $\Gamma\backslash\calH^2$ such that
\bqn
i(\mu,\delta_c)=\per(\gamma)
\eqn
whenever $c\subset \Gamma\backslash\calH^2$ is a closed geodesic and $\gamma$ a hyperbolic element representing $c$.
\end{thm}

Let $\ioo{x}{y}\subset\partial\calH^2$, respectively $\icc{x}{y}$, denote the open, respectively closed, interval with endpoints $x$ and $y$.
If $R=\ioo{d}{a}\times\ioo{b}{c}$ is an open rectangle in $\partial\calH^2$ defined by $(a,b,c,d)\in(\partial\calH^2)^{[4]}$

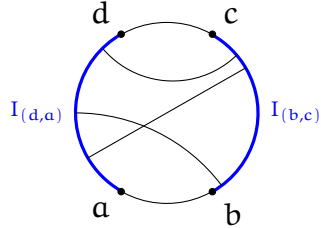
\begin{center}
\begin{tikzpicture}[scale=1.2]
\draw (0,0) circle [radius=1cm];
\draw[blue, very thick] (.5,-.87) arc[start angle=-60, end angle=60, radius=1];
\draw[blue, very thick] (-.5,.87) arc[start angle=120, end angle=240, radius=1];

\draw (-1,0) node [left,blue] {\tiny $\ioo{d}{a}$};
\draw (1,0) node [right,blue] {\tiny $\ioo{b}{c}$};


\filldraw (.5,.87) circle [radius=1pt] node[above right] {$c$};
\filldraw (.5,-.87) circle [radius=1pt] node[below right] {$b$};
\filldraw (-.5,.87) circle [radius=1pt] node[above left] {$d$};
\filldraw (-.5,-.87) circle [radius=1pt] node[below left] {$a$};

\draw (.87,.5) -- (-.87,-.5);
\draw (-.7,.71) arc [start angle=220, end angle=315, radius=1];
\draw (-1,0) arc [start angle=90, end angle=37, radius=2];

\end{tikzpicture}
\end{center}
we show in the course of the proof of Theorem~\ref{thm:per} that 
\bq\label{eq:mu1}
\mu(R)=\sup\left\{[a',b',c',d']:\,\icc{d'}{a'}\times\icc{b'}{c'}\subset R\text{ and }(a',b',c',d')\in\hgf\right\}\,.
\eq

\begin{remark}\label{rem:mu}  Notice that the properties \eqref{it:cr1} and \eqref{it:cr2} of the positive crossratio
reflect the fact respectively that $\mu$ is flip-invariant and that
if $x_1,\dots,x_5$ are positively oriented, the measure of the rectangle 
with vertices $x_1,x_2,x_4,x_5$ is related to the sum of the measures of the rectangle
with vertices $x_1,x_2,x_3,x_5$ and the one with vertices $x_1,x_3,x_4,x_5$.  

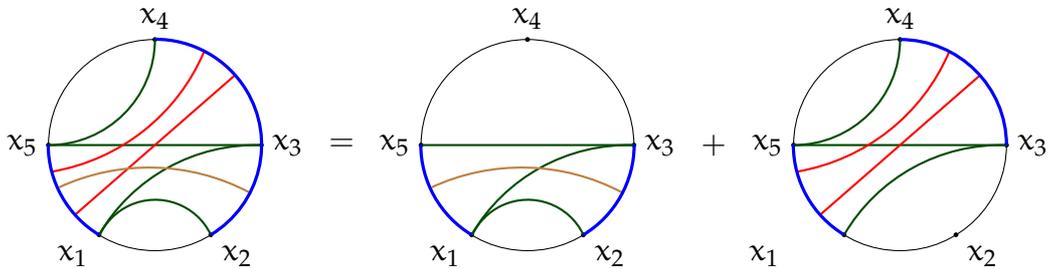
\begin{center}
\begin{tikzpicture}[scale=.7]
\draw (-7,0) circle [radius=2];
\draw (0,0) circle [radius=2];
\draw (7,0) circle [radius=2];

\draw (-3.5,0) node {$=$};
\draw (3.5,0) node {$+$};

\filldraw (-8.05, -1.7) circle [radius=1pt];
\draw (-8.05, -1.7) node[below left] {$x_1$};
\filldraw (-5.95, -1.7) circle [radius=1pt];
\draw (-5.95, -1.7) node[below right] {$x_2$};
\filldraw (-5, 0) circle [radius=1pt];
\draw (-5,0) node[right] {$x_3$};
\filldraw (-7,2) circle [radius=1pt];
\draw (-7,2) node[above] {$x_4$};
\filldraw (-9, 0) circle [radius=1pt];
\draw (-9,0) node[left] {$x_5$};

\draw[thick, black!70!green] (-9,0) arc [start angle=270, end angle=360, radius=2];
\draw[thick, black!70!green] (-9,0) -- (-5,0);
\draw[thick, black!70!green] (-5.95, -1.7) arc [start angle=25, end angle=155, radius=1.15];
\draw[thick, black!70!green] (-8.05, -1.7) arc [start angle=150, end angle=88, radius=3.4];
\draw[thick,red] (-8.5, -1.32) -- (-5.5,1.32);
\draw[thick, brown] (-8.8, -.8) arc [start angle=115, end angle=62, radius=4];
\draw[thick, red] (-8.92, -.5) arc [start angle=280, end angle=337, radius=3.8];

\draw[very thick, blue] (-8.05, -1.7) arc [start angle=238.5, end angle=180, radius=2];
\draw[very thick, blue] (-5.95, -1.7) arc [start angle=301.5, end angle=450, radius=2];

\filldraw (-1.05, -1.7) circle [radius=1pt];
\draw (-1.05, -1.7) node[below left] {$x_1$};
\filldraw (1.05, -1.7) circle [radius=1pt];
\draw (1.05, -1.7) node[below right] {$x_2$};
\filldraw (2, 0) circle [radius=1pt];
\draw (2,0) node[right] {$x_3$};
\filldraw (0,2) circle [radius=1pt];
\draw (0,2) node[above] {$x_4$};
\filldraw (-2, 0) circle [radius=1pt];
\draw (-2,0) node[left] {$x_5$};

\draw[thick, black!70!green] (-2,0) -- (2,0);
\draw[thick, black!70!green] (1.05, -1.7) arc [start angle=25, end angle=155, radius=1.15];
\draw[thick, black!70!green] (-1.05, -1.7) arc [start angle=150, end angle=88, radius=3.4];
\draw[thick, brown] (-1.8, -.8) arc [start angle=115, end angle=62, radius=4];

\draw[very thick, blue] (-1.05, -1.7) arc [start angle=238.5, end angle=180, radius=2];
\draw[very thick, blue] (1.05, -1.7) arc [start angle=301.5, end angle=360, radius=2];

\filldraw (5.95, -1.7) circle [radius=1pt];
\draw (4.945, -1.7) node[below left] {$x_1$};
\filldraw (8.05, -1.7) circle [radius=1pt];
\draw (8.05, -1.7) node[below right] {$x_2$};
\filldraw (9, 0) circle [radius=1pt];
\draw (9,0) node[right] {$x_3$};
\filldraw (7,2) circle [radius=1pt];
\draw (7,2) node[above] {$x_4$};
\filldraw (5,0) circle [radius=1pt];
\draw (5,0) node[left] {$x_5$};

\draw[thick, black!70!green] (5,0) arc [start angle=270, end angle=360, radius=2];
\draw[thick, black!70!green] (5,0) -- (9,0);
\draw[thick, black!70!green] (5.95, -1.7) arc [start angle=150, end angle=88, radius=3.4];
\draw[thick,red] (5.5, -1.32) -- (8.5,1.32);
\draw[thick, red] (5.08, -.5) arc [start angle=280, end angle=337, radius=3.8];

\draw[very thick, blue] (5.95, -1.7) arc [start angle=238.5, end angle=180, radius=2];
\draw[very thick, blue] (9,0) arc [start angle=0, end angle=90, radius=2];

 \end{tikzpicture}
 \end{center}

Notice however that there is no continuity assumption on a positive crossratio,
and as a result the resulting measure can have atoms,
 which is reflected by the following inequality valid for
  all  $(a,b,c,d)\in \hgf$:
 \bq\label{eq:mu2}
\mu(\ioo{d}{a}\times\ioo{b}{c})\leq[a,b,c,d]\leq\mu(\icc{d}{a}\times\icc{b}{c})\,.
\eq
As a matter of facts, \eqref{eq:mu1} and \eqref{eq:mu2} are equivalent characterizations of the measure $\mu$, \cite{BIPP_pos}.
\end{remark}

The positive crossratio defined above is obtained from an $\FF$-valued crossratio 
on pairwise transverse quadruples of Lagrangians $\calL(\ftn)^{\{4\}}$,
defined as follows:  if $(\ell_1,\ell_2,\ell_3,\ell_4)$ is a quadruple of $n$-dimensional vector subspaces in $\FF^{2n}$,
satisfying the properties that $\ell_1\pitchfork\ell_2$ and $\ell_3\pitchfork\ell_4$,
we define the crossratio
\bqn
\mathrm{CR}\colon\calL(\ftn)^{\{4\}}\to\FF\,,
\eqn
to be 
\bqn
CR(\ell_1,\ell_2,\ell_3,\ell_4)=\det(\mathrm{p}_{\ell_1}^{\parallel\ell_2}\circ\mathrm{p}_{\ell_3}^{\parallel\ell_4}|_{\ell_1})\,,
\eqn
where $\mathrm{p}_{\ell_i}^{\parallel\ell_j}$ is the projection on $\ell_i$ parallel to $\ell_j$.

Assume now that we are in the situation of Theorem~\ref{thm:repr_curr}.
In particular $\rho\colon\pis\to\Sp(2n,\FF)$ is a maximal representation.
According to Theorem~\ref{thm:max=fr} there is an equivariant map
\bqn
\varphi\colon \hs\to\calL(\ftn)
\eqn
sending positively oriented triples of points in $\hs$ to maximal triples of Lagrangians.
Then the formula 
\bqn
[x_0,x_1,x_2,x_3]_\rho:=-\nu(\mathrm{CR}(\varphi(x_1),\varphi(x_2),\varphi(x_3),\varphi(x_4)))
\eqn
defines a positive crossratio on $\hsf$ with values in $\nu(\FF^\times)$ 
and one then shows that the geodesic current $\mu_\rho$ given by Theorem~\ref{thm:per} satisfies the conclusion of Theorem~\ref{thm:repr_curr}.  

\subsection{Discrete valuations and atomic currents}

Since $\FF$ is real closed with valuation $\nu$, its value group $\nu(\FF^\times)$ is a $\QQ$-vector subspace of $\RR$.
Remarkably, if $\nu$ is non-Archimedean the crossratio $[\,\cdot\,,\,\cdot\,,\,\cdot\,,\,\cdot\,]_\rho$ takes its values in a smaller subgroup.

\begin{thm}[\cite{BIPP_pos}]\label{thm:val_cr}  Let $\rho\colon\pis\to\Sp(2n,\FF)$ be a maximal representation with $\FF$ non-Archimedean
and $v\colon\FF\to\RR\cup\{\infty\}$ a valuation.
Then there exists $K\in\NN_{\geq1}$, such that
\bqn
[x_1,x_2,x_3,x_4]_\rho\in\frac{1}{K}v(\QQ(\tr(\rho))^\times)
\eqn
for all $(x_1,x_2,x_3,x_4)\in H^{[4]}$.
\end{thm}

Together with the explicit relation between the crossratio and our construction of the geodesic current $\mu_\rho$ 
we obtain the following:

\begin{cor}\label{cor:val_cr}  Let $\rho\colon\pis\to\Sp(2n,\FF)$ be a maximal representation
and assume that $\QQ(\tr(\rho))$ has discrete valuation.  
Then the associated current $\mu_\rho$ is a {\em multicurve}, 
that is a finite sum of Dirac masses along closed geodesics and geodesics with endpoints in cusps.
\end{cor}

Indeed if $\QQ(\tr(\rho))$ has discrete valuation, we may assume that, up to scaling, the value group $\nu(\FF^\times)=\ZZ$
and hence, by Theorem~\ref{thm:val_cr}, the crossratio takes values in $\frac{1}{K}\ZZ$.
This implies by \eqref{eq:mu1} that the current is purely atomic with atoms whose mass is in $\frac{1}{K}\NN$.

\begin{ex}\label{ex:val} (Example \ref{ex:maxrep} continued)
  There are unique order compatible valuations $\nu_{0_+}$ and $\nu_\infty$ 
on $\QQ(X)$ such that 
\bqn
\nu_{0_+}\left(\frac{1}{X}\right)=-1\,,
\eqn
respectively 
\bqn
\nu_\infty(X)=-1\,.
\eqn
Thus these valuations extend uniquely to $\ov{\QQ(X)}^{0_+}$ and $\ov{\QQ(X)}^{\infty}$, respectively.
If $\rho$ is as in Example~\ref{ex:rho_closed}, since $\QQ(\tr(\rho))\subset\QQ(X)$, these valuations are discrete on $\QQ(\tr(\rho))$.
It follows then from Corollary~\ref{cor:val_cr} that the geodesic currents associated to $\rho_{0_+}$ and to $\rho_\infty$ are both multicurves.
\end{ex}

\subsection{The cases in which $\mu_\rho$ is a lamination}
\label{subsec:lam}
%
%

We give now a novel argument to deduce that if $n=1$ and $\partial\Sigma=\varnothing$, 
the length functions in the Thurston boundary are length functions of measured geodesic laminations,
that is of geodesic currents with vanishing self-intersection 
(see for example \cite[\S~8.1--8.3]{Martelli} or for example the original proof in \cite{Brum2}).
We deduce in fact directly from the properties of the crossratio the following:

\begin{prop}\label{prop:meas_lam} Let $\Sigma$ be a compact oriented surface
with possibly non-empty boundary.
Let $(\rho,\FF)$ represent a closed point in $\rsp\partial\Ximax(\Sigma,\SL(2,\RR))$.  
Then the associated geodesic current $\mu_\rho$ represents a measured geodesic lamination.
\end{prop}
 
 Indeed, if $(\ell_1,\ell_2,\ell_3,\ell_4)$  is a positively oriented quadruple in $\PP^1(\FF)$ and
 $x:=CR(\ell_1,\ell_2,\ell_3,\ell_4)$, then $x\geq1$ and $CR(\ell_4,\ell_1,\ell_2,\ell_3)=\frac{x}{x-1}\geq1$.
 Using that $\FF$ is non-Archimedean one sees that either $v(x)=0$ or $v\left(\frac{x}{x-1}\right)=0$.

 
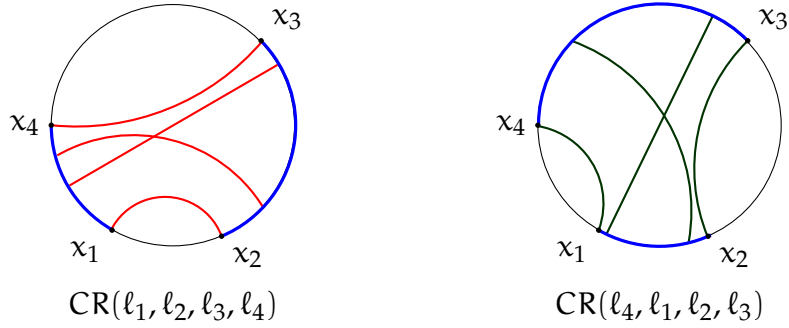
\begin{center}
 \begin{tikzpicture}[scale=.8]
 \draw (-4,0) circle [radius=2];
 \draw (4,0) circle [radius=2];

 \draw[thick, red] (-6,0) arc [start angle=264, end angle=320, radius=4];
 \draw[thick, red] (-5.7,-1) -- (-2.28,1);
 \draw[thick, red] (-5.9,-.5) arc [start angle=120, end angle=31, radius=2.5];
 \draw[thick, red] (.-3.2, -1.84) arc [start angle=20, end angle=155, radius=.99];

\draw[very thick, blue] (-6,0) arc [start angle=180, end angle=240, radius=2];
\draw[very thick, blue] (.-3.2, -1.84) arc [start angle=293, end angle=404, radius=2];

\draw[thick, black!80!green] (5.45, 1.4) arc [start angle=135, end angle=202, radius=3];
\draw[thick, black!80!green] (2.55, 1.4) arc [start angle=70, end angle=-10, radius=3];
\draw[thick, black!80!green] (3.13, -1.8) -- (4.87,1.8);
\draw[thick, black!80!green] (3, -1.74) arc [start angle=340, end angle=440, radius=1.3];

\draw[very thick, blue] (3, -1.74) arc [start angle=240, end angle=295, radius=2];
\draw[very thick, blue] (5.45, 1.4) arc [start angle=44, end angle=180, radius=2];

 \filldraw (-5, -1.74) circle [radius=1pt];
 \draw (-5, -1.74) node[below left] {$x_1$};
 \filldraw (-3.2, -1.84) circle [radius=1pt];
 \draw (.-3.2, -1.84) node[below right] {$x_2$};
 \filldraw (-2.55, 1.4) circle [radius=1pt];
 \draw (-2.55, 1.4) node[above right] {$x_3$};
 \filldraw (-6,0) circle [radius=1pt];
 \draw (-6,0) node[left] {$x_4$};
 
 \node at (-4,-3) {$CR(\ell_1,\ell_2,\ell_3,\ell_4)$};

 \filldraw (3, -1.74) circle [radius=1pt];
 \draw (3, -1.74) node[below left] {$x_1$}; 
 \filldraw (4.8, -1.84) circle [radius=1pt];
 \draw (4.8, -1.84) node[below right] {$x_2$};
 \filldraw (5.45, 1.4) circle [radius=1pt];
 \draw (5.45, 1.4) node[above right] {$x_3$};
 \filldraw (2,0) circle [radius=1pt];
 \draw (2,0) node[left] {$x_4$};
 
 \node at (4,-3) {$CR(\ell_4,\ell_1,\ell_2,\ell_3)$};

 \end{tikzpicture}
 \end{center}
  
 This implies by \eqref{eq:mu2}  that given any $(x_1,x_2,x_3,x_4)\in\hgf$, 
 either $\mu_\rho(\ioo{x_4}{x_1}\times\ioo{x_2}{x_3})=0$ or $\mu_\rho(\ioo{x_3}{x_4}\times\ioo{x_1}{x_2})=0$.
 Thus any two geodesics in the support of $\mu_\rho$ are either disjoint or coincide.

\medskip
Geodesic laminations also arise non-trivially in higher rank and,
when this is the case, we can find in the building
invariant  trees  that are convex for the distance
$\dnf$.
Let $\rho\colon\pis\to\Sp(2n,\FF)$ be a maximal representation defining a point in $\rsp\partial\Ximax(\Sigma,\Sp(2n,\RR))$;
in particular $\FF$ is non-Archimedean and admits a non-trivial valuation $\nu\colon\FF\to\RR$, 
with value group denoted by $\Lambda:=v(\FF^\times)$.
Let us assume that $(\rho,\FF)$ is a closed point 
and that the geodesic current $\mu_\rho$ associated to $\rho$ corresponds to a measured geodesic lamination,
which we regard as a $\pis$-invariant measured geodesic lamination of $\calH^2$.
Then on the set $\calT_{\mu_\rho}$ of complementary regions of this lamination
we have the distance function $d_{\mu_\rho}$ 
that assigns to two regions the $\mu_\rho$-measure of a geodesic segment connecting two interior points thereof.
Then $d_{\mu_\rho}$ is a $\Lambda$-valued distance function 
and $\calT_{\mu_\rho}$  is $0$-hyperbolic, that is embeddable in an
$\RR$-tree.
Under these conditions we have the following:

\begin{thm}[\cite{BIPP_pos}]\label{thm:last}  There is an equivariant embedding
\bqn
\calT_{\mu_\rho}\hookrightarrow\bsp
\eqn
of the $0$-hyperbolic $\Lambda$-metric space $\calT_{\mu_\rho}$ 
into the building $\bsp$, which is isometric with respect to $d_{\mu_\rho}$
and the distance $\dnf$ in \eqref{eq:dnf}.

This embedding extends to an isometric equivariant embedding of the $\RR$-tree dual to $\mu_\rho$
(obtained if there are no atoms by completion of the $\Lambda$-tree $\calT_ {\mu_\rho}$) in  the completion $\obsp$.
\end{thm}

\subsection{From the closed points in the real spectrum compactification to projectivized geodesic currents}\label{subsec:rsp_to_pr_curr}
If $(\rho,\FF)$ represents a point in $\rspclmax(\Sigma,\Sp(2n,\RR))$, 
according to Theorem~\ref{thm:repr_curr} there is an associated
geodesic current, once one fixes the valuation on $\FF$. 
 While if $\FF$ is Archimedean  we can -- and will --  take $\nu=-\ln$,  
when $\FF$ is  non-Archimedean the choice of the order compatible valuation is unique only up to positive multiples.
Moreover in light of Theorem~\ref{thm:bdry_cl} this current does not vanish.
Hence we obtain a well defined map
\bq\label{eq:rsp_to_pr_curr}
\ba
\rspclmax(\Sigma,\Sp(2n,\RR))&\stackrel{\calC}{\to}\PP(\calC(\mathring{\Sigma}))\\
(\rho,\FF)\,\,\quad\qquad&\longmapsto\,\,[\mu_\rho]
\ea
\eq
into the compact space of projectivized geodesic currents on
$\mathring{\Sigma}$.

Let $\mcgp$ be the orientation preserving mapping class group of $\Sigma$.
The group $\Out(\pis)$ acts on the character variety $\Xi(\pis,\Sp(2n,\RR))$ and 
$\mcgp$, seen as a subgroup of $\Out(\pis)$, leaves $\Ximax(\Sigma,\Sp(2n,\RR))$ invariant.
In addition $\mcgp$ also acts on the space $\calC(\mathring{\Sigma})$ of geodesic currents on $\mathring{\Sigma}$.

\begin{thm}[\cite{BIPP_max}]
  \label{thm:mcgp_action}
  The map in \eqref{eq:rsp_to_pr_curr}
\bqn
\rspclmax(\Sigma,\Sp(2n,\RR))\stackrel{\calC}{\to}\PP(\calC(\mathring{\Sigma}))
\eqn
is continuous and $\mcgp$-equivariant.
\end{thm}
This map is part of a natural diagram that involves the Weyl chamber length compactification.
Let $\thp\Ximax(\Sigma,\Sp(2n,\RR))$ be the closure of $\Ximax(\Sigma,\Sp(2n,\RR))$ in $\thp\Xi(\pis,\Sp(2n,\RR))$.
Identify $\fa$ with $\RR^n$ so that, under this identification,
${\fabp}$ corresponds to 
\bqn
\{x\in\RR^n:\,x_1\geq\dots\geq x_n\geq0\}\,.
\eqn
The norm on $\fa$ that corresponds to the multiplicative norm $N$ on $A$ defined in \eqref{eq:nf_ex}
is
\bqn
\|x\|_1:=\sum _{i=1}^n|x_i|
\eqn
and this induces a natural map 
\bqn
\xymatrix@1{
\PP({\fabp}^\pish)\ar[r]^{\|\,\cdot\,\|_1}
&\PP(\RR_{\geq0}^\pish)\,,
}
\eqn
where $\pish$ consists of the hyperbolic elements in $\pis$, inducing a map
\bqn
\xymatrix@1{
\thp\Ximax(\Sigma,\Sp(2n,\RR))\ar[r]^-{\|\,\cdot\,\|_1}
&\PP(\RR_{\geq0}^\pish)\,.
}
\eqn
The map in Theorem~\ref{thm:RSp->ThP} restricts to a continuous surjective map
\bq
\label{eq:mapRspClMaxToThP}
\xymatrix@1{
   \rspclmax(\Sigma,\Sp(2n,\RR))\ar[r]
&\,\thp\Ximax(\Sigma,\Sp(2n,\RR))
}
\eq
so that, using Theorem~\ref{thm:repr_curr} we obtain the following commutative diagram
\bqn
\xymatrix{
&\rspclmax(\Sigma,\Sp(2n,\RR))\ar[r]^-\calC\ar[r]\ar[dd]
&\PP(\calC(\mathring{\Sigma}))\ar[dd]^i
\\
\qquad\qquad\Ximax(\Sigma,\Sp(2n,\RR))\ar@{^{(}->}[ur]\ar@{^{(}->}[dr]
\\
&\thp\Ximax(\Sigma,\Sp(2n,\RR))\ar[r]^-{\|\,\cdot\,\|_1}
&\PP(\RR_{\geq0}^\pish)\,,
}
\eqn
where $i$ is the map induced by the intersection of a geodesic current in $\PP(\calC(\mathring{\Sigma}))$ 
with a closed geodesic associated to an element in $\pish$.

In the case of $n=1$ and $\partial\Sigma=\varnothing$, the image of $\calC$ is the compactification of Teichm\"uller space by
geodesic currents as described in Bonahon \cite{Bon88-curr}.

\subsection{Density of discrete points}
Observe that if $(\rho,\FF)$ is any representation and $\FF$ is of transcendence degree one over $\qbarr$,  
then $\QQ(\tr(\rho))$ has discrete valuation,
so that by Corollary~\ref{cor:val_cr}  the associated current $\mu_\rho$ is a {\em multicurve}.
 %
We obtain  a ``density of integral points" {\em \`a la Thurston}
combining this observation with the fact that 
the set of closed points in the real spectrum compatification of any semialgebraic subset of $\RR^n$
whose associated field has transcendence degree one over $\qbarr$ is dense, \cite[Proposition~4.2]{Brum2}:

\begin{cor}\label{cor:dense}
  The set of points in $\rspcl\partial\Ximax(\Sigma,\Sp(2n,\RR))$ 
for which the corresponding current is a multicurve is dense.
\end{cor}

This implies the analogous density of discrete points result  for
the Weyl chamber length compactification, using the continuous
surjective map \eqref{eq:mapRspClMaxToThP}.

\begin{cor}\label{cor:density in ThP}
  The set of points in $\thp\partial\Ximax(\Sigma,\Sp(2n,\RR))$ 
  for which the corresponding current is a multicurve is dense.
\end{cor}

\subsection{Towards amenability of the action of $\mcgp$ on the real spectrum compactification}\label{subsec:amen}
Observe that certain points in the Thurston boundary of Teichm\"uller space have large stabilizers in the mapping class group.
For instance the stabilizer of a measured geodesic lamination
supported on a separating closed geodesic 
contains the product of the mapping class groups of each component.  In contrast with this we have:

\begin{thm}[\cite{BIPP_max}]\label{thm:va_stab}  The action of $\mcgp$ on $\rsp\Ximax(\Sigma,\Sp(2n,\RR))$ has virtually Abelian stabilizers.
\end{thm}

The theorem is of course interesting only for the $\mcgp$-action on the boundary, 
since it acts properly discontinuously on $\Ximax(\Sigma,\Sp(2n,\RR))$.

The results in Theorem~\ref{thm:va_stab}  and the fact that Hamenst\"adt proved that the action of the $\mcgp$ on the compact metrizable Hausdorff space of complete measured geodesic laminations of $\Sigma$ is topologically amenable, except for some elementary surface types, 
\cite{Hamen_MCG}, lead to the following natural question:

\begin{q}  Does $\mcgp$  act topologically amenably on the compact metrizable space $\rspclmax(\Sigma,\Sp(2n,\RR))$?
\end{q}


\begin{proof}[Sketch of the proof of Theorem~\ref{thm:va_stab}]
The  very general property of the real spectrum of a ring  that every point in $\rsp\Ximax(\Sigma,\Sp(2n,\RR))$ contains a unique closed point in its closure
\cite[Proposition~7.1.24]{BCR} leads to a continuous retraction 
\bqn
\xymatrix@1{
\rsp\Ximax(\Sigma,\Sp(2n,\RR))\ar[r]
&\rspclmax(\Sigma,\Sp(2n,\RR))\,.
}
\eqn
Next, it $\Sigma_1\subset\Sigma$ is a subsurface, it follows from \cite{BIW} that 
restricting a homomorphism from $\pis$ to $\pi_1(\Sigma_1)$ induces a polynomial map
\bqn
\xymatrix@1{
\Ximax(\Sigma,\Sp(2n,\RR))\ar[r]
&\Ximax(\Sigma_1,\Sp(2n,\RR))\,.
}
\eqn
that extends continuously to a map
\bqn
\xymatrix@1{
\rsp\Ximax(\Sigma,\Sp(2n,\RR))\ar[r]
&\rsp\Ximax(\Sigma_1,\Sp(2n,\RR))
}
\eqn
that is equivariant with respect to any subgroup of $\mcgp$ that leaves $\Sigma_1$ invariant.
We caution the reader that however in general this map does not send closed points to closed points.

A subgroup $\Lambda<\mcgp$ stabilizing a point $[(\rho,\FF)]\in\rsp\Ximax(\Sigma,\Sp(2n,\RR))$,
will stabilize the corresponding closed point $[(\ov\rho,\ov\FF)]$ 
and hence will projectively stabilize the corresponding geodesic current $\mu_{\ov\rho}$ via the map $\calC$.
Applying the decomposition theorem \cite[Theorem~1.2\,(2)]{BIPP_curr},
there is a geodesic lamination $\calE\subset\mathring{\Sigma}$ consisting of simple closed geodesics
and such that the restriction of $\mu_{\ov\rho}$ to a complementary region $\calR\subset\mathring{\Sigma}\smallsetminus\calE$
either vanishes identically or its support fills $\calR$ in a specific way.
%

By passing if necessary to a subgroup $\Lambda'$ of $\Lambda$ leaving invariant all the components of $\mathring{\Sigma}\smallsetminus\calE$,
it would suffice to show that the projective stabilizer of the current restricted to each of the complementary regions is virtually Abelian. 
This is the case for the complementary regions of the latter kind, as the projective stabilizer of such current is either finite or cyclic infinite.
In the former case this is however not true, as the restriction of $\mu_{\ov\rho}$ to $\calR$ vanishes.
We consider then the subsurface with boundary $\ov\calR$ 
and use the restriction map
\bqn
\xymatrix@1{
\rsp\Ximax(\Sigma,\Sp(2n,\RR))\ar[r]
&\rsp\Ximax(\ov\calR,\Sp(2n,\RR))
}
\eqn
 applied to $[(\ov\rho,\ov\FF)]$ to produce an element $[(\pi,\KK)]\in\rsp\Ximax(\ov\calR,\Sp(2n,\RR))$  that is fixed by $\Lambda'$.
 Then the corresponding closed point $[(\ov\pi,\ov\KK)]\in\rspclmax(\ov\calR,\Sp(2n,\RR))$ is $\Lambda'$-invariant 
 and has an associated non-vanishing geodesic current.  
 In this way we reduce the complexity of the surface and the process stops after finitely many steps.
\end{proof}

\subsection
{A set of discontinuity for $\mcgp$ and harmonic maps}
%
Define the systole of a representation $\rho\colon\pis\to G_\FF$ by
\bqn
\Syst(\rho):=\inf_{\gamma\in\pish}
\ell_{N_\FF}(\rho(\gamma)) \,,
\eqn
where we recall that $\pish$ are the hyperbolic elements of $\pis$.
Observe that, in the real spectrum boundary of the character variety
of maximal representations, the positive systole subset 
\bqn
\rsp\Omega_\mathrm{max}(\Sigma,\Sp(2n,\RR)):=
\left\{(\rho,\FF)\in\rspcl\partial\Ximax(\Sigma,\Sp(2n,\RR)):\,\Syst(\rho)>0\right\}
\eqn
is well-defined and independent of the choice of the multiplicative norm $N\colon A\to\RR_{\geq1}$.
Using the map $\calC$ from the real spectrum to the projectivized geodesic
currents of Theorem~\ref{thm:mcgp_action} and the analogous result
for currents \cite{BIPP_curr}
we show:


\begin{cor}\label{cor:disc}  For $n>1$,
 the set $\rsp\Omega_\mathrm{max}(\Sigma,\Sp(2n,\RR))$ is  a non-empty open set of discontinuity
for the mapping class group action on the real spectrum
boundary $\rspcl\partial\Ximax(\Sigma,\Sp(2n,\RR))$.
\end{cor}
The fact that for $n>1$ the set is not empty follows from \cite[Corollary~1.10]{BIPP_curr} if $\partial\Sigma=\varnothing$;
to establish this fact if $\partial\Sigma\neq\varnothing$ we use the Strubel coordinates on $\Ximax(\Sigma,\Sp(2n,\RR))$, \cite{Strubel}.

\medskip
We now turn to an application to harmonic maps and, from this viewpoint, we will consider the CAT(0) distance
on the completion $\obsp$ induced by the embedding in Proposition~\ref{prop:isom_emb}.
As $\Ximax(\Sigma,\Sp(2n,\RR))$ is a closed semialgebraic set consisting of non-parabolic representations,
\cite{BIW}, Theorem~\ref{thm:proper} and Corollary~\ref{cor:harm} can be restated as follows:
\begin{cor}\label{cor:harmSp} Let $(\rho,\FF)\in\rspcl\partial\Ximax(\Sigma,\Sp(2n,\RR))$.  Then the $\pis$-action 
on the complete CAT(0) space $\obsp$ is KS-proper.
In particular if $\partial\Sigma=\varnothing$ for any hyperbolic metric $h$ on $\Sigma$
there exists a $\pis$-equivariant Lipschitz harmonic map 
\bq\label{eq:harm}
\xymatrix@1{
(\wt\Sigma,\wt h)\ar[r]
&\obsp\,.
}
\eq
\end{cor}
In fact, if $(\rho,\FF)$ represents a point in the set $\rsp\Omega_\mathrm{max}(\Sigma,\Sp(2n,\RR))$ of discontinuity of $\mcgp$ 
and $\partial\Sigma=\varnothing$ much more can be said: 

\begin{cor}  If $\partial\Sigma=\varnothing$, for every $(\rho,\FF)\in\rsp\Omega_\mathrm{max}(\Sigma,\Sp(2n,\RR)$,
there exists a hyperbolic metric $h$ on $\Sigma$ for which the $\pi(\Sigma)$-equivariant Lipschitz harmonic map in \eqref{eq:harm}
is also conformal.
\end{cor}

In fact, under the above conditions,
it follows from \cite[Corollary~1.5\,(2)]{BIPP_curr}
that the action of $\pis$ on $\obsp$ is displacing in the sense of \cite{DGLM}.
This implies by an argument of \cite{Lab08}, that the energy functional $E(h)$ on the Teichm\"uller space $\calT(\Sigma)$ of $\Sigma$ is proper.
Thus the map $\calT(\Sigma)\to\RR$ defined by $h\mapsto E(h)$ is proper and 
by \cite[Corollary~1.3]{Wentworth} it has a minimum which is conformal.

\begin{remark}  
A couple of spurious but interesting remarks related to the above considerations:
\be
\item Since the $\pis$-action on $\obsp$ is displacing, by \cite{DGLM} 
the $\pis$-orbit maps in $\obsp$ are quasi-isometric embeddings.
\item Our approach with working with semialgebraic sets that are defined by polynomial with coefficients in $\qbarr$
yields that the fields $\FF$ involved in the real spectrum are all countable.  
In particular the CAT(0) space $\obsp$  obtained as completion of a countable metric space is separable
and this might be useful to analyse regularity properties of harmonic maps.
\ee
\end{remark}

%

\section{Hitchin components and positive representations}\label{sec:Hpos}
In this section we want to describe briefly, in analogy with the case of maximal representations,
certain results concerning the real spectrum compactification of the Hitchin component $\XiHit(\pis,\PSL(n,\RR))$, 
where $\Sigma$ is a connected oriented compact surface with negative Euler characteristic;
in addition we assume that $\partial\Sigma=\varnothing$.
This restriction will be the main difference with our treatment of maximal representations.

\medskip
Recall that the Hitchin component $\XiHit(\pis,\PSL(n,\RR))$ is the quotient modulo the $\PSL(n,\RR)$-action 
of the connected component $\Hom_\mathrm{Hit}(\pis,\PSL(n,\RR))$ 
of $\Homred(\pis,\PSL(n,\RR))$
containing $\pi_n\circ\rho$, where $\rho\colon\pis\to\PSL(2,\RR)$ is any holonomy representation of a hyperbolic structure
and $\pi_n$ is the irreducible $n$-dimensional  representation of $\PSL(2,\RR)$.
Thus $\XiHit(\pis,\PSL(n,\RR))$ is a connected component of the $\PSL(n,\RR)$-character variety of $\pis$.
As all of these objects are semialgebraic defined by polynomials with $\qbarr$-coefficients,
we may consider their extension to any real closed field $\FF$
(see the discussion after Example~\ref{ex:orders} and \cite[Definition~5.1.2]{BCR}).
Since the operation of extension to $\FF$ sends semialgebraically connected components to
semialgebraically connected components, this unambigously leads to an $\FF$-connected component
$\Homhit(\pis,\PSL(n,\FF))\subset \Homred(\pis,\PSL(n,\FF))$, (see \cite[Proposition~5.3.6(ii)]{BCR}), 
whose quotient $\XiHit(\pis,\PSL(n,\FF))$ coincides with the corresponding extension 
$\XiHit(\pis,\PSL(n,\RR))_\FF$ of the semialgebraic set $\XiHit(\pis,\PSL(n,\RR))$.

It is possible to give a characterization of $\FF$-Hitchin representations,
that is representations in $\Hom_\mathrm{Hit}(\pis,\PSL(n,\FF))$, in terms of positivity.
Given any ordered field $\KK$, the concept of positive triple, resp. quadruple, of complete flags in $\KK^n$
can be defined in exactly the same way as for $\KK=\RR$, since it involves only positivity of triple-ratios, resp. crossratios.
As in \S~\ref{sec:max}, let $\hs\subset\partial\calH^2$ be the set of fixed points of non-trivial elements\footnote{Since $\partial\Sigma=\varnothing$, 
all elements in $\pis$ are hyperbolic.}
in $\pis<\PSL(2,\RR)$ and let $\calF(\KK^n)$ be the set of complete flags in $\KK^n$.

\begin{defi}  A representation $\rho\colon\pis\to\PSL(n,\KK)$ is {\em positive} 
if there exists an equivariant map $\varphi\colon\hs\to\calF(\KK^n)$ sending positively oriented
triples, resp. quadruples, of points in $\hs$ to positively oriented triples, resp. quadruples, of complete flags.
\end{defi}

In her thesis project  X.~Flamm is currently working out the analogue of Theorem~\ref{thm:max=fr}: 
a representation  $\rho\colon\pis\to\PSL(n,\FF)$ is an $\FF$-Hitchin representation if and only if it is positive.


The proof of such result uses the Bonahon--Dreyer coordinates, \cite{Bonahon_Dreyer} that give an explicit semialgebraic description 
of $\XiHit(\pis,\PSL(n,\RR))$, which extends to $\XiHit(\pis,\PSL(n,\FF))$ for any real closed field $\FF$.

Let $\rsp\XiHit(\pis,\PSL(n,\RR))$ be the closure of $\XiHit(\pis,\PSL(n,\RR))$ in the real spectrum compactification $\rsp\Xi(\pis,\PSL(n,\RR))$ and 
$\rspclhit(\pis,\PSL(n,\RR))$ the subset of closed points.  
We have then the following analogue of Theorem~\ref{thm:bdryMax}:

\begin{thm}\label{thm:3.3}
There is a bijective correspondence
between points in the real spectrum compactification $\rsp\XiHit(\pis,\PSL(n,\RR))$ and equivalence classes $(\rho,\FF)$,
where $\rho\colon\pis\to\PSL(n,\FF)$  is positive and $\FF$ is the real closure 
of the field $\QQ(\tr(\Ad(\rho)))$, where $\Ad$ is the adjoint representation of $\PSL(n,\FF)$.

In addition $(\rho,\FF)$ represents a closed point if and only if 
$\FF$ is Archimedean over the ring of traces $\QQ[\tr(\Ad(\rho))]$ of $\Ad\circ\rho$.
\end{thm}

\medskip
Turning to the symmetric space aspect, in analogy with the case of maximal representations,
we will directly describe the extension to $\FF$ of all relevant geometric objects.

Let $\calP^1(n)_\FF$ be as in \S~\ref{subsec:RSp}, 
\bqn
A_\FF:=\left\{D=\diag(\lambda_1,\dots,\lambda_n):\,\lambda_i\in\FF_{\geq0},\,\prod_{i=1}^n\lambda_i=1\right\}
\eqn
and let 
\bqn
\abp_\FF:=\{D\in A_\FF:\,\lambda_1\geq\dots\geq\lambda_n\}
\eqn
be the (multiplicative) closed Weyl chamber.
Then the formula in \eqref{eq:deltaF} provides us with an $\abp_\FF$-valued distance function $\delta_\FF$.

In contrast with the case of the symplectic group, we use the following multiplicative norm
\bqn
\ba
N_\FF\colon A_\FF&\longrightarrow\,\,\FF_{\geq1}\\
D&\longmapsto\max_{i\neq j}\frac{\lambda_i}{\lambda_j}\,.
\ea
\eqn
If $\FF$ admits a non-trivial order compatible valuation $v\colon\FF\to\RR\cup\{\infty\}$,
which -- we recall -- is the case for all fields occurring in the real spectrum of the character variety of $\pis$,
we obtain the distance (if $\FF\subset\RR$) or the pseudo-distance (if $\FF$ is non-Archimedean)
\bqn
d_{N_\FF}(A_1,A_2):=-v(N_\FF(\delta_\FF(A_1,A_2)))\,,
\eqn
for $A_1,A_2\in\calP^1(n)_\FF$.
If $\FF$ is non-Archimedean, we let $\bpsl$ be the quotient of $\calP^1(n)_\FF$ by the $d_{N_\FF}=0$ relation.
In all cases, that is both if $\FF$ is Archimedean or not, the translation length of an element $g\in\PSL(n,\FF)$ is given by
\bqn
\ell_{N_\FF}=-v\left(\frac{|\lambda_1|}{|\lambda_n|}\right)\,,
\eqn
where $\lambda_1,\dots,\lambda_n$ are the eigenvalues in $\FF[i]$ of a representative of $g$, 
ordered so that $|\lambda_1|\geq\dots\geq|\lambda_n|$.

Echoing the work of Martone and Zhang \cite{Martone_Zhang},  for $\FF=\RR$,
one can define an appropriate crossratio on $\calF(\FF^n)$ and use the abstract framework described in \S~\ref{sec:max} to show:

\begin{thm}\label{thm:4.4}  Let $\rho\colon\pis\to\PSL(,\FF)$ be a positive representation,
where $\FF$ is a real closed field endowed with a non-trivial valuation $v\colon\FF\to\RR\cup\{\infty\}$,
which we assume to be $-\ln$ if $\FF\subset\RR$.
Then there is a geodesic current $\mu_\rho$ on $\Sigma$ such that 
\bqn
i(\mu_\rho,\delta_c)=\ell_{N_\FF}(\rho(\gamma))
\eqn
for every closed geodesic $c\subset\Sigma$ and hyperbolic element $\gamma\in\pis$ representing $c$.
\end{thm}

The explicit relation between the crossratio and the corresponding geodesic current 
allows one to obtain the analogue of Theorem~\ref{thm:val_cr} and Corollary~\ref{cor:val_cr}. As a result we also get:

\begin{cor}\label{cor:density in ThP_Hit}
The set of points in $\thp\partial\XiHit(\pis,\Sp(2n,\RR))$ 
for which the corresponding current is a multicurve is dense.
\end{cor}
The above theorem thus gives us a map
\bqn
\xymatrix@1{
\rspclhit(\pis,\PSL(n,\RR))\ar[r]^-\calC
&\PP(\calC(\Sigma))
}
\eqn
into the projective space of geodesic currents that is $\mcgp$-equivariant
and which can be shown to be continuous.  As in the case of maximal representations, 
the set
\bqn
\rsp\Omega_\mathrm{Hit}(\pis,\PSL(n,\RR)):=
\left\{(\rho,\FF)\in\rspcl\partial\XiHit(\pis,\PSL(n,\RR)):\,\Syst(\rho)>0\right\}
\eqn
is an open set of discontinuity for the action of the mapping class group.
It follows from \cite{BIPP_curr} that:

\begin{cor}  For $n\geq1$, the set $\rsp\Omega_\mathrm{Hit}(\pis,\PSL(n,\RR))$ is not empty.
\end{cor}

Turning to the applications to harmonic maps, we consider the CAT(0) metric on the completion $\obpsl$
induced by the embedding in Proposition~\ref{prop:isom_emb}.  
In complete analogy with the case of maximal representations we have:

\begin{cor}\label{cor:4.6}  \be
\item  If $(\rho,\FF)\in\rspcl\partial\XiHit(\pis,\PSL(n,\RR))$ the $\pis$-action on the complete CAT(0) space
$\obpsl$ is KS-proper.  
In particular for any hyperbolic metric $h$ on $\Sigma$ there exists a $\pis$-equivariant Lipschitz harmonic map
\bqn
(\wt\Sigma,\wt h)\longrightarrow\obpsl\,.
\eqn
\item For every $(\rho,\FF)\in\rsp\Omega_\mathrm{Hit}(\pis,\PSL(n,\RR))$ there exists a hyperbolic metric $h$ on $\Sigma$
for which the $\pis$-equivariant Lipschitz harmonic map in (1) is also conformal.
\ee
\end{cor}

\begin{remark} Here is a couple of loose ends:
\be
\item The analogue of Theorem~\ref{thm:last}, namely the question whether in case $\mu_\rho$ is a geodesic currents 
its associated tree $\calT_{\mu_\rho}$ embeds into $\bpsl$ presents new difficulties and is open in general.
\item The fact that the action of $\mcgp$ has elementary stabilizers, that is the analogue of Theorem~\ref{thm:smallstab} and of Theorem~\ref{thm:va_stab},
require a flexible enough notion of positivity for representations of $\pis$ with $\partial\Sigma\neq\varnothing$
and is a current topic of investigation.
\ee
\end{remark}

\bibliographystyle{alpha}
\bibliography{refs}

\end{document}